\documentclass[a4paper,10pt]{article}%[eqsection]
\everymath{\displaystyle}
\usepackage{amssymb,amsmath,amsthm}
\usepackage{graphicx}

\usepackage{amsfonts}
\usepackage{hyperref}
\usepackage{latexsym}
\usepackage{amsthm}
\usepackage{amsmath}
\usepackage{amssymb}
\usepackage{amsopn}
\usepackage{geometry}
\usepackage{amscd}
\usepackage{mathrsfs}
\usepackage{graphicx}
\usepackage[english]{babel}
\usepackage[english]{babel}
\newcommand{\Om}{\Omega}

\newenvironment{pf}{\noindent{\sc Proof}.\enspace}{\rule{2mm}{2mm}\medskip}
\newenvironment{pfn}{\noindent{\sc Proof} \enspace}{\rule{2mm}{2mm}\medskip}
\newtheorem{theorem}{Theorem}[section]
\newtheorem{proposition}{Proposition}[section]
\newtheorem{lemma}{Lemma}[section]
\newtheorem{corollary}{Corollary}[section]

\newtheorem{remark}{Remark}[section]
\newtheorem{remarks}{Remark}[section]
\newtheorem{definition}{Definition}[section]
\newcommand{\be}{\begin{equation}}
\newcommand{\ee}{\end{equation}}
\newcommand{\teta}{\theta}
\newcommand{\om}{\omega}

\newcommand{\e}{\varepsilon}

\newcommand{\ov}{\overline}
\newcommand{\wtilde}{\widetilde}
\newcommand{\R}{\mathbb R}
\newcommand{\C}{\mathbb C}

\newcommand{\Z}{\mathbb Z}

\newcommand{\N}{\mathbb N}
\newcommand{\T}{\mathbb T}

\renewcommand{\a }{\alpha }
\renewcommand{\b }{\beta }
\newcommand{\s }{\sigma }
\newcommand{\ii }{{\rm i} }
\renewcommand{\d }{\delta }
\newcommand{\D }{\Delta}
\newcommand{\g }{\gamma}

\renewcommand{\l }{\lambda }

\newcommand{\vphi}{\varphi }

\renewcommand{\t }{\tau }
\renewcommand{\o }{\omega }

\newcommand{\norma}{|\!\!|}

\newcommand{\lin}{{\cal L}}
\newcommand{\matr}{{\cal M}}

\newcommand{\nors}[1]{|\!\!| #1  |\!\!|_{s}}
\newcommand{\norsone}[1]{|\!\!| #1  |\!\!|_{s_1}}

\begin{document}

\title{{\bf Sobolev quasi periodic solutions of multidimensional  wave equations
% NLW on $ \T^d $
%of 
with a multiplicative potential
%nonlinear wave equations  in any dimension
}}

\date{}

\author{Massimiliano Berti, Philippe Bolle}

\maketitle

\noindent
{\bf Abstract:}
We prove the existence of  
quasi-periodic solutions  for wave equations with a multiplicative potential 
on $ \T^d $, $ d \geq 1 $, and finitely differentiable nonlinearities, quasi-periodically forced in time. 
The only external parameter is the length of the frequency vector. The solutions have Sobolev regularity both
in time and space. The proof is based on a Nash-Moser  iterative scheme as in \cite{BB10}. 
% scheme developed for  NLS 
The key % as in \cite{BB10} % under the weakest
 tame estimates for the inverse linearized operators 
 are obtained by a multiscale inductive argument, 
 % as in \cite{BB10}. % for NLS
%  along scales os Sobolev spaces.  The proof 
which  is more difficult than for  NLS due to  the dispersion relation of the wave equation. 
We prove the ``separation properties" of the small divisors 
assuming weaker non-resonance conditions than in  \cite{B5}.
\\[2mm]
{\it Keywords:} Nonlinear wave equation, 
Nash-Moser Theory, KAM for PDE, Quasi-Periodic Solutions, Small Divisors, 
Infinite Dimensional Hamiltonian Systems.
%\footnote{Supported by MIUR 
%``Variational Methods and Nonlinear Differential Equations" and by the European Research Council under FP7.}.
\\[1mm]
2000AMS subject classification: 35Q55, 37K55, 37K50.

\section{Introduction}

The first  existence results  of quasi-periodic solutions for Hamiltonian PDE
were proved by  Kuksin \cite{K1} and Wayne \cite{Wa1} 
for one dimensional  (1-$d$) nonlinear wave (NLW) and 
nonlinear Schr\"odinger (NLS) equations, extending 
 KAM theory. This approach consists  in generating iteratively 
a sequence of canonical changes of variables which bring the Hamiltonian into a normal form  with an 
invariant torus at the origin.
This  procedure requires, at each step, to invert linear ``homological equations", 
which have  constant coefficients %($=$ reducible)
and can be  solved by  imposing the ``second order Melnikov" non-resonance conditions.
The final KAM torus is linearly stable.
These pioneering results  were limited to 
Dirichlet boundary conditions because 
the eigenvalues of $ \partial_{xx} $ 
 had to be simple:
the second order Melnikov non resonance conditions are 
violated already for periodic boundary conditions. 

In such a case, the first existence results of quasi-periodic solutions 
were proved by Bourgain \cite{Bo1} 
extending the  approach of Craig-Wayne  \cite{CW} for periodic solutions. 
The search of the embedded torus is reduced to solving 
a functional equation in scales of Banach spaces,  by some 
Newton implicit function procedure. 
The main advantage of this scheme is to require only the ``first order Melnikov"
non-resonance conditions  to solve the homological equations.
These conditions are essentially  the minimal non-resonance assumptions.
Translated in the KAM language 
 this  corresponds to allow a  normal form with non-constant coefficients around the torus.
The main difficulty is that 
the homological equations are  PDEs with {\it non-constant} coefficients and 
are small perturbations of a diagonal operator having 
arbitrarily small eigenvalues. 

\smallskip

At present, the theory for $1$-$d$ NLS and NLW equations 
has been sufficiently understood
(see e.g. \cite{K2}, \cite{Po2},  \cite{KP}, \cite{Po3}, \cite{C}, \cite{BBi})
but much work remains 
in higher space dimensions.
The main difficulties are:
\begin{enumerate}
\item the eigenvalues  
 of  $ - \Delta + V(x) $ appear in clusters of unbounded sizes,
\item the eigenfunctions are ``not localized with respect to the exponentials".
\end{enumerate}

Roughly speaking, an eigenfunction $ \psi_j $ of $ - \Delta + V(x) $ is localized 
with respect to the exponentials, if its Fourier coefficients $ (\hat \psi_{j})_i $ rapidly  converge  to zero
(when $ |i - j | \to  \infty $).  
This property  always holds in $ 1 $ space dimension (see \cite{CW}) but 
may fail for $ d \geq 2 $, see \cite{B3}. It implies that the matrix which represents 
(in the eigenfunctions basis) the multiplication
operator for an analytic function has an exponentially fast decay off the diagonal. 
%From a dynamical point of view 
It reflects into a ``weak interaction" between  different 
``clusters of small divisors".
Problem 2  has been often bypassed 
replacing the multiplicative potential $ V(x) $ by a ``convolution potential"  
$ V * ( e^{\ii j \cdot x}) := m_j  e^{\ii j \cdot x}$,  $ m_j \in \R $, $ j \in \Z^d $. 
The ``Fourier multipliers"  $ m_j $   play the role of ``external parameters".

\smallskip

The first existence results of quasi-periodic solutions for analytic NLS and NLW  like
\be\label{NLSNLW}
\frac{1}{\ii} u_t = B u + \e \partial_{\bar u} H(u, \bar u) \, , \quad 
u_{tt} + B^2 u + \e F'(u) = 0 \, , \quad x \in \T^d \, , \ \  d \geq 2 \, , 
\ee
where $ B $ is a Fourier multiplier,  have been proved by Bourgain  
 \cite{B3}, \cite{B5}, by extending the Newton approach in \cite{Bo1} (see also \cite{B4} for periodic solutions). 
Actually this scheme is very convenient to overcome problem 1, because 
it requires  only the  first order Melnikov non-resonance conditions and therefore does not exclude multiplicity of
normal frequencies (eigenvalues). The main difficulty concerns the multiscale inductive argument to 
estimate the off diagonal exponential decay of the inverse linearized operators 
 in presence of huge clusters of small divisors. The proof is based on 
a repeated use of the resolvent identity and  
fine techniques of subharmonicity and semi-algebraic set theory, essentially to obtain refined 
 measure and  ``complexity" estimates for sublevels of functions. 

\smallskip

Also the KAM approach has been recently extended by  Eliasson-Kuksin \cite{EK} 
for  NLS on $ \T^d $ with Fourier multipliers and analytic nonlinearities. The key  issue is to 
control more accurately the perturbed frequencies after the KAM iteration and, in this way,
verify the  second order Melnikov non-resonance conditions, we refer 
also to \cite{GXY}, \cite{PX}, \cite{BBP1} for related techniques.  
We also mention  \cite{EK1} which proves the reducibility of a linear
Schr\"odinger equation % externally 
 forced by a small  multiplicative potential,  quasi-periodic in time.

On the other hand, a similar reducibility KAM result for NLW on $ \T^d $ is still an  open 
problem: the possibility of imposing the second order Melnikov conditions for wave equations in higher space dimensions  is
still uncertain. 

\smallskip

In the recent paper  \cite{BB10} we proved 
the existence of quasi-periodic solutions  for quasi-periodically forced NLS on $ \T^d $ 
with  finitely differentiable nonlinearities 
 (all the previous results were valid for analytic nonlinearities, actually polynomials in \cite{B3}, \cite{B5})
and a multiplicative potential $ V(x) $ (not small). Clearly  a difficulty 
is that the matrix which represents the multiplication operator 
has only a polynomial decay off the diagonal, and not exponential. 
The proof is based on a Nash-Moser  iterative scheme  in Sobolev scales
(developed for periodic solutions also in \cite{BB}, \cite{BB07},  \cite{BBP}, \cite{BP}) 
and novel techniques  for  estimating the high Sobolev norms of the solutions of the (non-constant coefficients) homological equations.
In particular we assumed that $ - \Delta + V(x) > 0 $ in order to prove the ``measure and complexity" estimates
by means of elementary eigenvalue variations arguments,  avoiding 
subharmonicity and semi-algebraic techniques as in \cite{B5}.  

\smallskip

The goal of this paper is to prove an analogous result -see Theorem  \ref{thm:main}-
for $ d $-dimensional nonlinear  wave equations with a quasiperiodic-in-time nonlinearity
like
\be\label{eq:main} 
u_{tt} - \Delta u + V(x) u  = \e f(\om t, x,u)   \, , \quad x \in {\T}^d  \, ,  \, \ \e >  0 \, , 
\ee 
where the multiplicative potential  $ V$ is in $ C^{q} ({\T}^d;\R) $, 
 $ \om \in \R^\nu $ is a non-resonant frequency vector (see \eqref{diophan0}, \eqref{diofgr}),   
% the nonlinearity is  quasi-periodic in time 
and 
\be\label{nonli}
f \in C^q ( {\T}^\nu \times {\T}^d \times \R;\R)   
\ee
for some  $ q \in \N $ large enough (fixed in Theorem \ref{thm:main}). 
The NLW equation is more difficult 
than  NLS because the small divisors 
stay near a  cone,  see \eqref{def:Y}, 
and not a paraboloid. 
Therefore it is harder to prove the ``separation properties" of the Fourier indices of the small divisors, see 
section \ref{sec:sepa}.  
In this paper we use 
a non-resonance condition which is weaker than in Bourgain \cite{B5}, see remark \ref{rem:bou}.
After the statement of Theorem \ref{thm:main} 
we explain the other main differences with respect to  \cite{B5} and
 \cite{BB10}.

\smallskip

Concerning the potential we suppose that 
\be\label{nonke}
{\rm Ker} (- \Delta + V(x)) = 0 \, . 
\ee
\begin{remark}
In \cite{BB10} we assumed the stronger condition  $- \Delta + V(x) > 0 $.  See comments after 
${\it Theorem} \, \ref{thm:main}$. Note that also  in \eqref{NLSNLW}  the Fourier operator $ B^2 > 0 $  is  positive. 
\end{remark}
In \eqref{eq:main} we use only one external parameter, namely the length of the frequency vector 
(time scaling). More precisely 
we assume that the frequency vector $ \om $ is co-linear  with a fixed 
 vector $ \bar \o \in \R^{\nu}  $, 
\be\label{baromega}
\om = \l \bar \om \, , \quad \l \in \Lambda := [1/2, 3/2] \subset \R \, ,  \quad |\bar \om | \leq 1 \, , 
\ee
where $ \bar \omega $ is Diophantine, namely  for some $  \g_0 \in (0,1) $, 
\be\label{diophan0}
|\bar \om \cdot l | \geq \frac{2 \g_0}{ | l |^{\nu}} \, , \quad \forall l \in \Z^{\nu} \setminus \{ 0 \} \, , 
\ee
and 
\be\label{diofgr}
\Big| \sum_{1 \leq i\leq j \leq \nu} {\bar \om}_i {\bar \om}_j p_{ij} \Big|  \geq \frac{\g_0}{ | p |^{\t_0}} \, , 
\quad \forall p \in  \Z^{\frac{\nu ( \nu + 1 )}{2}} \setminus \{ 0 \}
\, . 
\ee
There exists $ \bar \o $ satisfying \eqref{diophan0} and \eqref{diofgr} at least for  $ \t_0 >  \nu (\nu+1)   -1 $ and $ \g_0 $ small, see Lemma \ref{diofp}.
For definiteness we fix $ \t_0 :=  \nu (\nu+1)  $.
\begin{remark}
For  NLS equations \cite{BB10} only condition \eqref{diophan0} is required, see comments after Theorem \ref{thm:main}.
\end{remark}

 \smallskip

The dynamics of the linear wave equation
\be\label{eq:L}
u_{tt} - \Delta u + V(x) u = 0 
\ee
is well understood.  The eigenfunctions of 
$$
(- \Delta + V(x)) \psi_j(x) = \mu_j \psi_j(x) 
$$
form a Hilbert basis in $ L^2 (\T^d )$ and  the eigenvalues  $ \mu_j \to + \infty $ as $ j \to + \infty $.
By  assumption \eqref{nonke} all the eigenvalues $ \mu_j$ are different from $0$.
We list them in non-decreasing order 
\be\label{eigDelta}
\mu_1 \leq \ldots \leq \mu_{n^-}  < 0 < \mu_{n^- +1} \leq \ldots 
\ee
where 
$ n^- $ denotes the number of negative eigenvalues (counted with multiplicity).

All the solutions of \eqref{eq:L} are the linear superpositions
of normal mode oscillations, namely 
$$
u(t,x) = \sum_{j=1}^{n^-} 
(\b_j^- e^{-\sqrt{|\mu_j|} t} + \b_j^+ e^{\sqrt{|\mu_j|} t}) \psi_j (x) +
\sum_{j \geq n^-+1} {\rm Re} ( a_j e^{\ii \sqrt{\mu_j} t} ) \, \psi_j (x) \, , \ \b_j^\pm  \in \R \, , a_j \in \C \, .
$$
The first $ n^- $ eigenfunctions correspond to hyperbolic directions where the dynamics is attractive/repulsive.
The other  infinitely many  eigenfunctions correspond to elliptic directions.

%On the other  infinitely many  eigenfunctions
%the solutions are periodic, quasi-periodic or almost periodic in time. 

\begin{itemize}
\item {\sc Question:} for $ \e $ small enough, do there exist quasi-periodic solutions of the nonlinear wave equation 
(\ref{eq:main}) for positive measure sets of 
$ \l \in [1/2, 3/2] $?
\end{itemize}

\noindent
Note that, if $ f(\vphi,x,0) \not\equiv 0 $ then $ u = 0 $ is not a solution of (\ref{eq:main}) for $ \e \neq  0 $.

\smallskip

The above question amounts to look for $(2 \pi)^{d+ \nu} $-periodic solutions $ u(\vphi,x) $ of
\be\label{eq:freq} 
(\om \cdot \partial_\vphi)^2 u  - \Delta u + V(x) u = \e f( \vphi , x, u)  
\ee 
in the Sobolev space
\begin{eqnarray}\label{def:Hs}
H^s  :=  H^s ( \T^\nu \times \T^d;  {\R}) &:=&  \Big\{  u(\vphi,x) := \sum_{ (l,j) \in \Z^\nu \times \Z^d}  u_{l,j} 
e^{\ii (l \cdot \vphi + j \cdot x)}  :  \, \| u \|_s^2 := K_0 \sum_{i \in \Z^{\nu + d}}  |u_{i}|^2 \langle i \rangle^{2s} < + \infty \, , 
\nonumber \\
& & \ \ \, u_{-i}=\ov{u_i} \  , \   {\rm where} \  \ i := (l,j)  \, ,  \  \langle i \rangle := \max(|l|,|j|,1) 
\Big\} 
\end{eqnarray}
for some $  (\nu + d) \slash 2 < s \leq q $. Above $ |j| := \max \{|j_1|, \ldots , |j_d| \}$.   For the sequel we fix 
$ s_0 > (d+\nu )\slash 2 $ so that  
$ H^s (\T^{\nu+d} ) \hookrightarrow  L^\infty (\T^{\nu+d} ) $, $ \forall s \geq s_0 $. 
The constant $ K_0 >  0  $ in  \eqref{def:Hs}  
is fixed (large enough) so that  $ |u|_{L^\infty}  \leq \| u \|_{s_0}$ and the  interpolation inequality 
\be \label{interp}
\| u_1 u_2 \|_s \leq \frac12 \| u_1 \|_{s_0} \| u_2 \|_s + \frac{C(s)}{2}  \| u_1 \|_s \| u_2 \|_{s_0}   \, , \quad 
\forall s \geq s_0 \, , \ u_1, u_2 \in H^s \, , 
\ee
holds with $ C(s) \geq 1 $, $ \forall s \geq s_0 $, and $ C(s) = 1 $, $ \forall s \in [s_0, s_1 ] $; the constant $ s_1 := s_1 (d, \nu) $ 
is defined in \eqref{Sgr}. 

\smallskip

The main result of the paper is:

\begin{theorem}  \label{thm:main} Assume \eqref{diophan0}-\eqref{diofgr}. 
There are $ s := s(d, \nu) $,  $ q := q(d, \nu) \in \N $, such that:
$ \forall f  \in C^q $, $ \forall V \in C^q $ satisfying \eqref{nonke}, 
 $ \forall \e \in [0, \e_0) $ small enough, there is a  map 
\be\label{normd}
u(\e, \cdot) \in C^1(\Lambda ;H^s) 
\quad {\rm with }  \quad \sup_{\l \in \Lambda} \| u (\e, \l) \|_s \to 0  \  {\rm as} \  \e \to 0 \,, 
\ee
and a Cantor like set $ {\cal C}_\e \subset  \Lambda := [1/2, 3/2] $ of 
asymptotically full Lebesgue measure, i.e.   
\be\label{Cmeas}
| {\cal C}_\e |  \to 1 \quad {\rm as } \quad \e \to 0  ,  
\ee
such that, $ \forall \l \in {\cal C}_\e $, $ u (\e, \l)  $ is a solution of (\ref{eq:freq}) with $ \om = \l \bar \om $.
Moreover, if $ V, f $ are of class $ C^\infty $ then $\forall \l\, , \ u(\e,\l) \in C^\infty (\T^d \times \T^\nu; \R)$. 
\end{theorem}

\noindent
Let us make some  comments on the result. 

\begin{enumerate}
\item 
The main novelties of Theorem \ref{thm:main} with respect 
to previous literature (i.e. \cite{B5}) are that we prove
 the existence of quasi-periodic solutions  for quasi-periodically forced NLW on $ \T^d $, $ d \geq 2 $, with a
\begin{description}
\item (i) \,  {\it multiplicative} finitely differentiable potential $ V(x) $, 
\item (ii) \, finitely {\it differentiable} nonlinearity, see \eqref{nonli}, 
\item (iii)  {\it pre-assigned} direction of the tangential frequencies, see \eqref{baromega}.
\end{description}
Moreover we weaken the non-resonance assumptions to ensure the separation properties of the small divisors.   
Theorem \ref{thm:main}  generalizes \cite{BB}
to the case of quasi-periodic solutions.
% and  also \cite{PY}, \cite{Fo} for wave equations in higher dimension. 
\item 
We underline that the present Nash-Moser approach 
requires essentially no  information about the localization of the eigenfunctions of  $ - \Delta + V(x) $
which, on the contrary, seem to be unavoidable to prove also reducibility with a KAM scheme. 
Along the multiscale analysis we use (as in \cite{BB10}) the exponential basis which diagonalizes $ - \Delta + m $ where $ m $
is the average of $ V(x) $. The key 
is to define  ``very regular" sites, namely take the constant  $ \Theta $
in Definition \ref{regulars}  large enough, 
depending on the potential $ V(x) $. 
In this way the number of sites to be considered as ``singular" increases. However, 
the separation properties of the  ``singular"  sites obtained in Lemma \ref{Bourgain} hold for any $ \Theta > 0 $, 
and this is sufficient for the applicability of the present multiscale approach. 
% are so  weak that they still hold, see section \ref{sec:sepa}. 
\item
Throughout this paper $ \e \in [0, \e_0] $ is fixed (small) and 
$ \l \in [1/2,3/2] $ is the only external parameter in equation \eqref{eq:main}. 
Then the bound \eqref{Cmeas} is an improvement with respect to the analogous Theorem 1.1  in \cite{BB10} 
(for NLS) where we only proved the existence of quasi-periodic solutions for a  Cantor set,
with asymptotically full measure, in the parameters $ (\e,\l ) \in [0,\e_0) \times [1/2, 3/2 ] $. 
%(this improvement  could be clearly obtained also for the NLS). 
\item 
We have not tried to optimize 
the estimates for $ q := q(d, \nu) $ and 
$ s := s(d, \nu) $.  In \cite{BB07} we proved the existence of periodic solutions in 
$ H^s_t H^1_x $ with $ s > 1 /2 $,  for  one dimensional NLW  equations
with  nonlinearities of class $ C^6 $, see the bounds (1.9), (4.28) in \cite{BB07}.
\end{enumerate}

\noindent
Let us make some comments about the proof.
The main differences with  respect to \cite{BB10} and \cite{B5} are:

\begin{enumerate}
\item Since we do not assume  that
$ - \Delta + V(x) $ is positive definite  (as in \cite{BB10}), but only the weaker assumption
\eqref{nonke}, the measure and complexity arguments in section \ref{sec:measure} are more 
difficult than in \cite{BB10}, section 6. The main reason why we can allow a finite
number of negative eigenvalues $ \mu_j < 0 $ in \eqref{eigDelta} is that the corresponding small divisors 
satisfy 
$$ 
- (\om \cdot l )^2 + \mu_j \leq \mu_j \leq \mu_{n^-} < 0 \, , \ \  \forall l \in \Z^\nu, \,  j = 1, \ldots, n^-, 
$$ 
namely are {\sc not} small, it is used in Lemma \ref{A-1}. 
The positivity of $ - \Delta + V(x)  $ was used in \cite{BB10} to prove the
measure and complexity estimates. Assuming only \eqref{nonke},  
the main difference concerns Lemma \ref{cor2} that we tackle with a Lyapunov-Schmidt type argument.
Note that Lemma \ref{cor2} only holds for  $ j_0 \notin {\cal Q}_N $ defined in \eqref{calQN}  
(in such a case the spectrum of the restricted operator $ \Pi_{N,j_0} ( - \D + V(x))_{E_{N,j_0}} $ in
\eqref{PNpr} is far away from zero by Lemma \ref{lzero}). 
This fact  requires to modify also the definition of $ N $-good sites, see 
Definition \ref{GBsite}, with respect to the analogous Definition 5.1 of \cite{BB10}. 
\item The separation properties of the small divisors  in section \ref{sec:sepa} are proved
under the non-resonance assumption $ {\bf (NR)} $ (see \eqref{NRom}, \eqref{diofgr}), 
%which is weaker  with respect to that in Bourgain \cite{B5}, Chapter 20:  
which is  a Diophantine condition 
for polynomials in $ \om $ of degree $ 2 $, while the condition in \cite{B5} for polynomials of higher degree, 
see remark \ref{rem:bou}. 
A Diophantine condition like  ${\bf (NR)} $  is necessary 
because % we are proving for This is due to the dispersion relation indeed we have to prove 
the singular sites are  integer points near a cone, see \eqref{singusites}, and not a paraboloid like for NLS. 
Then it is necessary to assume an irrationality condition on the ``slopes" of this cone.
Assumption $ {\bf (NR)} $ % in  \eqrefom0} 
seems to be the weakest possible. 
The improvement is in  the % different 
proof of  Lemma \ref{Bourgain} (different with respect to Lemma 20.14 of Bourgain \cite{B5}) 
% In Lemma \ref{Bourgain} 
which extends, to the quasi-periodic case, the arguments of \cite{BB}.
\item
Another technical simplification of the present approach with respect to \cite{B5}, Chapter 20, 
is to study NLW in configuration space without regarding \eqref{eq:main}  as a first order 
Hamiltonian complex 
system. The main difficulty concerns the measure estimates:  
the derivative with respect to $ \teta $ of the matrix in \eqref{matrpar} is not positive definite 
(this affects Lemmata \ref{lem:complexity} and, especially, \ref{cor2}). 
The main technical trick that we use 
is  the change of variables  \eqref{changevaria}.  We mention that  also 
Bourgain-Wang \cite{BW1}, section 6, deals with NLW in configuration space, where
the measure and complexity estimates are verified using  subharmonicity and 
semi-algebraic techniques. 
\end{enumerate}

\noindent
{\bf Acknowledgments:}  We thank Luca Biasco and Pietro Baldi for useful comments. 

\section{The linearized equation}\setcounter{equation}{0}

We  look for solutions of the NLW equation (\ref{eq:freq}) in $ H^s $ 
by means of a Nash-Moser iterative scheme.
The main step concerns the invertibility  of (any finite dimensional restriction of) the linearized operator 
\be\label{Linve}
{\cal L} (u) := {\cal L}(\o, \e,  u) := L_\om - \e g (\vphi,x)
\ee
where 
\be\label{pq}
L_{\om} :=  (\o \cdot  \partial_\vphi)^2 -   \Delta + V(x) \qquad {\rm and} \qquad  
g(\vphi,x) :=  (\partial_u f)(\vphi,x, u) \, .
\ee
%The function $ g $ depends also on $ \e, \l $ through $ u $. 
We decompose the multiplicative potential as
$$
V(x) = m + V_0 (x)
$$
where $ m $ is the average of $ V(x) $ and $ V_0 (x) $ has zero mean value. Then we write 
\be\label{Lomega} 
L_\om = D_\om +  V_0(x) \qquad {\rm where} \qquad D_\om :=  (\o \cdot  \partial_\vphi)^2 -   \Delta + m 
\ee
has  constant coefficients. % in the Fourier basis $(e^{\ii(l \cdot \varphi + j \cdot x)})$.  
In the Fourier basis $(e^{\ii(l \cdot \varphi + j \cdot x)})$,  
the  operator  ${\cal L}(u) $ 
is represented by the infinite dimensional self-adjoint matrix 
$$
A(\om) := A ( \om, \e,  u) := D  +  T   
$$
where 
$$
% D := {\rm diag} \d_i   := - (\l \bar{\om} \cdot l)^2 + \|j\|^2  + m   \, . 
 D :=  {\rm diag}_{(l, j) \in \Z^\nu \times \Z^d} \,  - (\om \cdot l)^2 + \|j\|^2 + m := {\rm diag}_{i \in \Z^b} \d_i  \, , 
 %  \, , \quad \d_i :=  - (\om \cdot l)^2 + \|j\|^2 + m \, , 
$$ 
\be\label{def:euc}
\| j \|^2 := j_1^2+ \ldots + j_d^2 \, , \quad   i := (l,j) \in \Z^{b}  := \Z^\nu \times \Z^d  \, , \quad \d_i := - (\om \cdot l)^2 + \|j\|^2 + m
\ee
and 
\be \label{Tmatrix}
T := %T(\e,\l) :=
T_2 - \e T_1  \, , \ \ 
T := (T_i^{i'})_{i ,  i'\in \Z^b} \, , \ \
T_i^{i'} :=  (V_0)_{j-j'}  - \e g_{i-i'}  
\ee 
represents the % {\it T\"oplitz} matrix which represents the 
multiplication operator by $ V_0(x) - \e g(\vphi,x) $.
The matrix $ T $ is  {\it T\"oplitz}, namely $ T_i^{i'} $ depends only on the difference of the indices $ i - i' $, and, 
since the functions $ g $, $ V \in H^s $,   then $ T_i^{i'} \to 0  $ as $ |i- i'| \to \infty $  at a polynomial rate.  

\smallskip

Along the iterative scheme of section \ref{sec:NM}, the function $ u $ (hence $g$) will depend on $(\e,\l)$, so that 
$T:=T(\e,\l)$ will be considered as a family of operators 
(or of infinite dimensional matrices representing them in the Fourier basis)
parametrized by $(\e,\l)$. Introducing an additional parameter $\theta$, we  
 consider the  family of infinite dimensional matrices 
\be \label{matrpar}
A(\e,\l,\theta)=D(\teta) +  T(\e,\l) 
\ee 
where
\be\label{def:Y}
D (\theta) :=  D(\l,  \teta) := {\rm diag}_{i \in \Z^b} \Big( - (\l  \bar \om \cdot l + \theta )^2 + \|j\|^2 + m \Big)  
\ee
and  $\norsone{T} + \norsone{\partial_{\lambda} T} \leq C $, depending on $ V $
(the norm $ \norsone{\ } $ is introduced in Definition \ref{defnormatr}). 
The main goal of the following sections is to prove polynomial off-diagonal decay for the inverse
of  the $ (2N+1)^b $-dimensional sub-matrices of $ A(\e, \l, \theta) $ centered at $ (l_0, j_0) $ denoted by
\be\label{ANl0}
A_{N,l_0, j_0}(\e, \l, \theta) := A_{|l - l_0| \leq N, |j - j_0| \leq N}(\e, \l, \theta) 
\ee
where 
$ |l| :=  \max \{|l_1|, \ldots, |l_\nu|\} $, $ |j| := \max \{|j_1|, \ldots, |j_d|\} $.  
The relation with $ \| j \| $ defined in \eqref{def:euc} is 
\be\label{supeuc}
|j| \leq \|j\| \leq \sqrt{d} |j| \, . 
\ee
If $ l_0 = 0 $ we use the simpler notation
$$
A_{N, j_0}(\e, \l, \theta) :=  A_{N,0,j_0} (\e, \l, \theta) \, .
$$
If also $ j_0 = 0 $,  we simply write
$$
A_{N}(\e, \l, \theta) :=  A_{N,0} (\e, \l, \theta) \, ,
$$
and, for $ \teta = 0 $, we denote
$$
A_{N,j_0}(\e, \l) :=  A_{N,j_0} (\e, \l, 0) \, .  
$$
By  \eqref{ANl0},  \eqref{matrpar}, \eqref{def:Y} and since $ T $ is T\"oplitz, 
the following {\it crucial} covariance property (exploited in Lemma \ref{Ntime}) holds:
\be\label{shifted}
A_{N, l_1,  j_1 } (\e, \l, \theta) = A_{N, j_1} (\e, \l , \teta +  \l  \bar{\om} \cdot l_1 ) \, . 
\ee

\subsection{Matrices with off-diagonal decay}\label{sec:off}
 
For  $ B \subset \Z^b  $ we introduce the subspace
$$
H^s_B := \Big\{ u = \sum_{i \in \Z^b}  u_{i} e_i  
\in H^s  \, : \, u_{i} = 0 \ {\rm if} 
\ i \notin B \Big\} 
$$
where  $ e_i  := e^{\ii (l \cdot \vphi + j \cdot x)}  $. When $ B $ is finite, the space $ H^s_B $ does not 
depend on $ s $ and will be denoted $ H_B $.   
For $ B, C \subset \Z^b  $ finite, 
we identify the space $ \lin^B_C $ of the linear maps $ L : H_B \to H_C $  with the space of matrices
$$
\matr^B_C := \Big\{ M = (M^{i'}_i)_{i' \in B, i \in C} \, , \ M^{i'}_{i} \in \C  \Big\} 
$$ 
identifying $ L $ with the matrix $ M $ with entries
$ M_i^{i'}  :=  (L e_{i'}, e_i)_0 $ % $ i' \in B $, $ i \in C $ (
where $ ( \, , \, )_0 := (2\pi)^{-b} ( \, , \, )_{L^2}$ denotes the normalized $ L^2 $-scalar product.
We consider also the $ L^2 $-operatorial norm 
\be\label{L2norm}
\| M^B_C \|_0 := \sup_{h \in H_B, h \neq 0} \frac{\| M^B_C h \|_0}{\| h \|_0}  \, .
\ee

\begin{definition} \label{defnormatr} {\bf ($s$-norm)}
The $ s $-norm of a matrix $ M \in \matr^B_C $ is defined by
$$
\norma M \norma_s^2 := K_0 \sum_{n \in \Z^b} [M(n)]^2 \langle n \rangle^{2s}  
$$
where $ \langle n \rangle := \max (|n|,1)$ (see \eqref{def:Hs}), % $ |n| := \sup_{p = 1, \ldots, b} |n_p| $, 
$$
[M(n)] := \begin{cases}
\max_{i-i'=n} |M^{i'}_i|  \ \ \quad   {\rm if}  \ \   n \in  C- B\\
0  \qquad \qquad \qquad  {\rm if}  \ \  n \notin C - B 
\end{cases}
$$
and the  constant $ K_0 > 0 $ is the one of \eqref{def:Hs}.
\end{definition}

The $ s $-norm is modeled on matrices  which represent the multiplication operator.

\begin{lemma}\label{lem:multi}
The (T\"oplitz) matrix $ T $ which represents the multiplication 
operator by %$ T_2 $ in \eqref{pq}, \eqref{Lomega} 
$ g  \in H^s $ satisfies
$\nors{T} \leq C \| g \|_s $.  % \quad \nors{T_2} \leq  \| V \|_s  $
\end{lemma}

In analogy with the  operators of multiplication by a function, 
the matrices with finite $ s $-norm satisfy interpolation inequalities (see \cite{BB10}). 
As a particular case, we can  derive from (\ref{interp})

\begin{lemma}\label{sobonorm} {\bf (Sobolev norm)}
$ \forall s \geq s_0 $ there is $ C(s) \geq 1 $ such that, for any finite subset $ B, C \subset \Z^b $, 
\be\label{opernorm}
\| Mw \|_s \leq (1/2) \norma M\norma_{s_0} \|w\|_s + (C(s)/2)  \norma M \norma_s \|w\|_{s_0} \, , \quad 
\forall M \in \matr^B_C \, , \ w \in H_B \, .
\ee
\end{lemma}

\subsection{A spectral lemma}

We denote % the space of (functions of the $ x $-variable only)
\be\label{E0}
E_{N,j_0} := \Big\{ u(x) := \sum_{|j-j_0| \leq N}  u_j e^{\ii j \cdot x} \, , \ u_j \in \C  \Big\} 
\ee
 (functions of the $ x $-variable only) and the corresponding orthogonal projector 
\be\label{proNj0}
\Pi_{N,j_0} : H^{s_0}(\T^d) \to E_{N,j_0} \, .
\ee
More generally, for a finite non empty subset $ B \subset \Z^d$ we denote by $\Pi_B$ the $L^2$-orthogonal projector onto the space 
$ E_B \subset L^2 (\T^d ) $ spanned by $\{ e^{\ii j \cdot x} \, : \, j \in B \}$.

We now prove a result on the spectrum of the restricted self-adjoint operator
\be\label{DVB}
(-\Delta + V)_B := \Pi_B (-\Delta + V)_{|E_B} \, 
\ee 
that shall be used for the measure estimates of Lemma \ref{cor2}.

We shall denote  (with a slight abuse of notation)
$$
\partial B := \Big\{ j \in B \ : \ {\rm d}(j, \Z^d \backslash B ) =1 \Big\} 
$$
where $ {\rm d}(j, j')  := | j - j' | $ denotes the distance associated to the sup-norm.
Note that,  if 
$ {\rm d}(0, \partial B) \geq L_0 $,  $ L_0 \in \N $,  % (note that $ {\rm d}(0, \partial B) \in \N $) 
then: either  
$$ 
{\mathtt B}(0, L_0-1)  := \{ j \in \Z^d \, : \, |j| \leq L_0 - 1 \}  \subset \Z^d \backslash B 
\qquad
{\rm or} \qquad
{\mathtt B}(0, L_0)  \subset  B 	\, . 
$$
%$ {\rm d}(0, \partial B) < L_0 $ (note that $ {\rm d}(0, \partial B) \in \N $) implies that 
%$$ 
%{\mathtt B}(0, L_0)  := \{ j \in \Z^d \, : \, |j| \leq L_0  \} \not \subset  B \qquad
%{\rm and} \qquad
%{\mathtt B}(0, L_0-1) \not \subset \Z^d \backslash B \, .
%$$
Recall \eqref{eigDelta} where 
$ n^- $ is the number of negative eigenvalues of $ -\Delta +V(x) $ (counted with multiplicity).

\begin{lemma} \label{lzero}
Let  
$ \beta_0 := \min \{ |\mu_{n^-}| / 2 , \mu_{n^- +1} \} $.
There is $ L_0 \in \N $, such that, if $ {\rm d}(0, \partial B) \geq L_0 $, then 
\begin{enumerate}
\item
if $ {\mathtt B}(0, L_0-1)  \subset \Z^d \setminus B $, then $(-\Delta + V)_B \geq \beta_0 I $,
\item
if $ {\mathtt B}(0, L_0) \subset B $, then $ (-\Delta + V)_B$ has $n^-$ negative eigenvalues, all of them $\leq -\beta_0$.
All the other eigenvalues of $(-\Delta + V)_B$ are $ \geq \beta_0 $. 
\end{enumerate}
\end{lemma}

\begin{pfn}
The  eigenvalues   \eqref{eigDelta} of $ - \Delta + V $ satisfy the min-max characterization
\be \label{minmax1}
\mu_p=\inf_{ G \subset H^1(\T^d) , \atop \dim G=p} \sup_{u \in G, \|u\|_{L^2}=1} Q(u) \, , \quad p = 1, 2, \ldots 
\ee
where $ Q : H^1 (\T^d; \R ) \to \R $ is the quadratic form 
\be\label{Quadrat}
Q(u) := \| \nabla u \|_{L^2}^2 + \int_{\T^d} V(x) u^2 (x) d x  
\ee
and the infimum in \eqref{minmax1} is taken over the subspaces $ G $ of $H^1(\T^d) $ 
of dimension $ p $.

Let $ {\cal H}^- \subset H^1 (\T^d) $ be the $ n^- $-dimensional orthogonal sum of the  eigenspaces associated to the negative 
eigenvalues $ \mu_1, \ldots , \mu_{n^-}$. Then
$$
Q(u) \leq \mu_{n^-} \|u\|^2_{L^2} % \stackrel{\eqref{beta0}} 
\leq -2 \beta_0 \|u\|^2_{L^2} \, , \quad \forall u \in {\cal H}^-  \, , 
$$
by the definition of $ \b_0 $. 
Moreover there is $ L_1 $ (large) such that $ G^- := \Pi_{L_1,0} {\cal H}^-$ (recall \eqref{proNj0}) has  dimension $n^-$ and 
\be\label{negatif}
Q(u)  \leq -\beta_0 \|u\|^2_{L^2} \, , \quad \forall u \in G^- \,  .
\ee
Let
\be\label{L0}
L_0 := 
\max \{ L_1, (\beta_0 + |V|_{L^\infty})^{1/2}\} \, . % ,\quad 
%|V|_{L^\infty} := \max_{x \in \T^d} |V(x)| \, . 
\ee
$1$) Assume ${\mathtt B}(0,L_0-1) \subset \Z^d \setminus B $. Then (using that $ {\rm d}(0,  B) \geq L_0 $)
$$
\|\nabla u \|^2_{L^2}  % = \sum_{j \in B} |u_j|^2 \| j \|^2
\geq L_0^2 \|u \|^2_{L^2} \,  , \quad \forall u \in E_B \, ,
$$
and, by \eqref{Quadrat},
$$
Q(u) \geq (L_0^2 -  |V|_{L^\infty} ) \| u \|^2_{L^2} 
\stackrel{\eqref{L0}} \geq \beta_0 \| u \|^2_{L^2} \, , \quad \forall u \in E_B \,  .
$$
Hence $ (-\Delta + V)_B \geq \beta_0 I$. 
\\[1mm]
$2$) Assume ${\mathtt B}(0,L_0) \subset B $. Let $ (\mu_{B,p}) $ be the non-decreasing sequence of the
eigenvalues of the self-adjoint operator $ (-\Delta + V)_B $, counted with multiplicity. 
They satisfy a variational characterization analogous to \eqref{minmax1} with the only difference 
that the infimum  is taken over the subspaces $ G \subset E_B  $. % of dimension $ p $.
Since $ {\mathtt B}(0,L_1) \subset {\mathtt B}(0, L_0) \subset B $, the subspace
$ G^- \subset E_B $ and, recalling that $ \dim G^- = n^- $,
$$
\mu_{B,n^-}= \inf_{G \subset E_B, \atop {\rm dim} G = n^-} \sup_{u \in G, \|u\|_{L^2}=1} Q(u) \leq 
\sup_{u \in G^-, \|u \|_{L^2}=1} Q(u) \stackrel{\eqref{negatif}} \leq -\beta_0 \, . 
$$
Moreover 
\begin{eqnarray*}
\mu_{B,n^-+1} & = & \inf_{G \subset E_B , \atop {\rm dim} G=n^-+1 } \sup_{u \in G, \| u \|_{L^2}=1} Q(u) \\
& \geq & \inf_{G \subset H^1(\T^d), \atop {\rm dim} G = n^- + 1} \sup_{u \in G, \| u \|_{L^2}=1} Q(u) \stackrel{\eqref{minmax1}} =  
\mu_{n^-+1}   \geq \beta_0 
\end{eqnarray*}
by the definition of $ \b_0 $. The proof of the lemma is complete.  
\end{pfn}

\section{The multiscale analysis} 
\setcounter{equation}{0}\label{multiscale}

We recall the multiscale Proposition \ref{propinv} proved in \cite{BB10}.
Given $ \Omega, \Omega' \subset  E \subset \Z^b  $ we define 
$$
{\rm diam}(E) := \sup_{i,i' \in E} |i-i'| \, , \qquad 
 {\rm d}(\Omega, \Omega') := \inf_{i \in \Omega, i' \in \Omega'} |i-i'|  \, . 
$$
Let $ \d \in (0,1)$ be fixed. 
\begin{definition}\label{goodmatrix}
{\bf ($N$-good/bad matrix)} The matrix $ A \in {\cal M}_E^E $, with $ E \subset \Z^b $, 
$ {\rm diam}(E) \leq 4 N $,  is $ N $-good if $ A $ is invertible and
$$
\forall s \in [s_0, s_1] \   , \ \ \nors{A^{-1}} \leq N^{\tau'+\d s}.
$$
Otherwise $ A $ is $ N $-bad.
\end{definition}

\begin{definition}\label{regulars} {\bf (Regular/Singular site)} Fix $ \Theta  \geq 1  $.
The index $ i \in \Z^{b}  $  is  {\sc regular} for $A=A(\e , \l , \theta)$ 
if $ |A_i^i| \geq \Theta $. Otherwise $ i $ is  {\sc singular}.
\end{definition}

\begin{definition}\label{ANreg}
{\bf ($(A,N)$-good/bad site)}
For $ A \in \matr^E_E $, we say that $ i \in E \subset \Z^b   $ is
\begin{itemize}
\item $(A,N)$-{\sc regular} if there is $ F \subset E$ such that
${\rm diam}(F) \leq 4N$,  ${\rm d}(i, E\backslash F) \geq N/2$ and
$A_F^F$ is $N$-good. 
\item $(A,N)$-{\sc good}  if it is regular for $A$ or $(A,N)$-regular. Otherwise we say that $ i $ is $(A,N)$-{\sc bad}.
\end{itemize}
\end{definition}

Let us consider  the new larger scale 
\be\label{newscale}
N' = N^\chi 
\ee
with $ \chi > 1 $.  For a matrix $ A \in \matr_E^E  $ we define $ {\rm Diag}(A) := ( \d_{ii'} A_i^{i'})_{i, i' \in E} $.

\begin{proposition} {\bf (Multiscale step, see  {\bf \cite{BB10})}} \label{propinv}
Assume
\be\label{dtC}
\d \in (0,1/2) \, ,  \ \tau' > 2 \tau + b + 1 \, , \ C_1 \geq 2 \, , 
\ee
and, setting $ \kappa := \tau' + b + s_0 $,
\be\label{chi1}
\chi (\t' - 2 \t - b) >  3 (\kappa + (s_0+ b) C_1 ) \, , \ 
\chi \delta > C_1 \, ,
\ee
\be\label{s1}
S \geq  s_1 > 3 \kappa + \chi (\tau + b) + C_1 s_0 \, . 
\ee
 $ \Upsilon > 0 $ being fixed , there exists $  N_0 (\Upsilon, S) \in \N $, $ \Theta (\Upsilon, s_1) > 0 $ large enough
(see Definition \ref{regulars}),  such that:
\\[1mm]
$ \forall N \geq N_0(\Upsilon, S) $, 
$ \forall E \subset \Z^b  $ with  ${\rm diam}(E) \leq 4N'=4N^\chi $, if $ A \in \matr_E^E $ satisfies 
\begin{itemize}
\item 
{\bf (H1)} $\norsone{A- {\rm Diag}(A)} \leq \Upsilon $ 
\item 
{\bf (H2)} $ \| A^{-1} \|_0 \leq (N')^{\tau}$
\item
{\bf (H3)} 
There is a partition of the  $(A,N)$-bad sites $ B = \cup_{\alpha} \Omega_\alpha$ with
\be\label{sepabad}
{\rm diam}(\Omega_\alpha) \leq N^{C_1} \, , \quad {\rm d}(\Omega_\alpha , \Omega_\beta) \geq N^2 \ , \ \forall \alpha \neq \beta \, ,
\ee
\end{itemize}
then $ A $ is $ N' $-good. More precisely
$$
\forall s \in [s_0,S]  \ , \  \  \nors{A^{-1}} \leq \frac{1}{4} ({N'})^{\tau' } \Big( ({N'})^{\d s}+ \nors{A- {\rm Diag}(A)} \Big) \,. 
$$
\end{proposition}

We shall apply Proposition \ref{propinv} to 
finite dimensional matrices  $ A_{N,i_0} $ (recall the notation in (\ref{ANl0})) 
which are 
obtained as restrictions of the infinite dimensional matrix $ A (\e,\l,\teta ) $ in (\ref{matrpar}).
It is convenient to introduce a notion of $ N $-good site for an infinite dimensional matrix. 

Let
\be\label{calQN}
{\cal Q}_N := \Big\{ j \in \Z^d \, : \, {\rm d}(0, \partial (j + [-N,N]^d)) < L_0 \Big\} \, , \ \ 
\check {\cal Q}_N := \Big\{ i=(l,j) \in \Z^d \, : \, j \in {\cal Q}_N  \Big\} 
\ee
where $L_0$ is defined in Lemma \ref{lzero}. We  shall always assume that $ N - 2 L_0 \geq N/2  $. 
\begin{definition} \label{GBsite} {\bf ($N$-good/bad site)}
A site $ i \in \Z^b $ is: 
\begin{itemize}
\item $ N$-{\sc regular } if $ A_{N,i} $ is $ N $-good (Definition \ref{goodmatrix}). Otherwise we say that $ i $ is $ N$-{\sc singular}.
\item $N$-{\sc good} if   $ i $  is  regular  (Definition \ref{regulars})
  {\rm or}  \hbox{for all} $ M \in \{ N- 2 L_0, N \} $, 
 all   the  sites   $ i' $  with  $ |i' - i| \leq M $ and  $ i' \notin \check {\cal Q}_M $  are $ M $-{regular}. 
 Otherwise, we say  that  $ i $ is $ N $-{\sc bad}. 
\end{itemize}
\end{definition}

Definition \ref{GBsite} is designed in view of the application of Proposition \ref{propinv}, because
we have

\begin{lemma} \label{ANgood}
Let $A=A_{N',i_0}$ with $i_0 \notin \check {\cal Q}_{N'}$. Then any $N$-good site 
$i \in i_0+[-N',N']^{d+\nu}$ is $(A,N)$-good. 
\end{lemma}

\begin{pf}
We decompose
\be\label{inizio}
E := i_0 + [-N',N']^{\nu + d} = G \times H  \quad  {\rm where} \quad
G := \Pi_{p=1}^\nu [a_p, b_p]  \, ,  \ H := \Pi_{q=1}^d [c_q,d_q]
\ee
and, writing  $ i_0 = (l_0, j_0 ) $,
$$
a_p := (l_0)_p-N'  \, , \ b_p :=(l_0)_p+N' \, , \ c_q :=(j_0)_q-N' \, , \ d_q := (j_0)_q+N' \, .
$$
Consider any $ N $-good site $  i := (l,j) \in  E $ (see Definition \ref{GBsite}).
If $i$ is a regular site, there is nothing to prove.  If $i$ is singular, 
we introduce its neighborhood 
\be\label{neiF}
F_N := F_N(i) :=G_N \times H_N \subset E \quad{\rm where} \quad
G_N := \Pi_{p=1}^\nu I_p   \subset G \, , \quad
H_N :=\Pi_{q=1}^d J_q   \subset H \, , \quad 
\ee
and the intervals  $ I_p \subset [a_p,b_p] $, $ J_q \subset [c_q,d_q] $ % \subset \Z $  
are defined as follows:
\begin{itemize}
\item if $l_p-a_p>N$ and $b_p-l_p>N$ (resp.
$j_q-c_q>N$ and $d_q-j_q>N$), then $I_p :=[l_p-N, l_p+N]$ (resp. $J_q :=[j_q-N, j_q +N]$);
\item if $l_p-a_p \leq N $ (resp. $j_q-c_q \leq N$), then 
$I_p :=[a_p, a_p+2N]$ (resp. $J_q :=[c_q, c_q+2N]$);
\item if $b_p-l_p\leq  N$ ( resp. $d_q-j_q \leq N$), then
$I_p :=[b_p-2N,b_p]$ (resp. $J_q :=[d_q-2N, d_q]$). 
\end{itemize}
By construction we have
\be\label{EFd}
{\rm d}(i, E \setminus F_N) \geq N 
\ee
and we can write % By the definition of $ F_N $ we can write  
\be\label{ritra}
F_N = \bar{\imath} + [-N,N]^{\nu +d} \quad {\rm for \ some } \ \ 
\bar{ \imath}= (\bar{l}, \bar{\jmath})  \in E \ \ {\rm with}  
\ \ |i - \bar{\imath}| \leq N  \, . 
\ee
For $ M = N - 2 L_0 $, we define as in \eqref{neiF} 
the sets  $ F_M := G_M \times H_M $, $ G_M :=  \Pi_{p=1}^\nu I_{M,p} $, $ H_M := \Pi_{q=1}^d J_{M,q} $,
and we write 
\be\label{itM}
F_M = \tilde{\imath} + [-M,M]^{\nu +d}   \quad {\rm for \ some } \ \ 
\tilde{\imath}=(\tilde{l}, \tilde{\jmath}) \ \ {\rm with} \ \ |i - \tilde{\imath}| \leq M  \, .  
\ee
We claim that
\be \label{dd}
{\rm d}( \partial H_N \backslash \partial H , H_M) \geq 2L_0 \, .
\ee
In fact, assume $j' \in  \partial H_N \backslash \partial H $. Then 
there is some $ q \in \{ 1, \ldots, d \} $ such that $ j'_q \in \partial J_q \backslash  \{ c_q, d_q \}$. 
By construction,  it is easy to see that $ {\rm d}(J_{M,q}, [c_q,d_q] \backslash J_q) \geq  2L_0+1$. 
Hence $ {\rm d}(j'_q, J_{M,q}) \geq 2L_0$ and ${\rm d}(j',H_M) \geq 2 L_0 $, proving \eqref{dd}.
 
We are now in position to prove that $ i $ is $ ( A, N )$-good. We distinguish two cases: 
\begin{itemize}
\item[(i)]
 $Ê{\rm d}(0, \partial H_N) \geq L_0 $. Since $ H_N = \bar \jmath + [-N,N]^d $ (see \eqref{neiF}-\eqref{ritra}) we get
 $ \bar \jmath \notin {\cal Q}_N $ (see \eqref{calQN}), namely  
 $\bar{\imath} \notin \check {\cal Q}_N$. Since $ i $ is a singular $N$-good site (see Definition \ref{GBsite}),  
 $|i - \bar{\imath}| \leq N$ (see \eqref{ritra}),  $\bar{\imath} \notin \check {\cal Q}_N$, 
 we deduce  that the matrix $A_{N, \bar{\imath}} = A_{F_N}^{F_N}$ is $ N $-good . 
As a consequence, since $ F_N \subset E $ (see \eqref{neiF}),  $ {\rm diam}(F_N) = 2 N $ (see \eqref{ritra})
$ {\rm d}(i, E \backslash F_N ) \geq N $ (see \eqref{EFd}), 
the site $ i $ is $(A,N)$-good  (see Definition  \ref{ANreg}). 
\item[(ii)] ${\rm d}(0, \partial H_N) < L_0$. 
It is an assumption of the Lemma that $ i_0 = (l_0, j_0) \notin \check {\cal Q}_{N'} $ which 
means ${\rm d}(0, \partial H) \geq L_0 $ (by \eqref{inizio} we have $ H = j_0 + [-N',N']^d $). Hence 
 ${\rm d}(0, \partial H_N \backslash \partial H) = {\rm d}(0, \partial H_N)  < L_0 $.
Hence, by (\ref{dd}), we deduce $ {\rm d}(0, H_M) \geq L_0 $ and therefore
$ {\rm d}(0, \partial H_M) \geq L_0 $. Then $\tilde{\imath} \notin \check {\cal Q}_M$ (the 
site $ \tilde{\imath}  $ is defined in \eqref{itM} and we have $ H_M = \tilde \jmath + [-M,M]^d $). 
Since $ i $ is singular and $N$-good, $ |i - \tilde{\imath}| \leq M $ (see \eqref{itM}), $\tilde{\imath} \notin \check {\cal Q}_M $, 
then the matrix $A_{M, \tilde{\imath}} = A_{F_M}^{F_M}$ is $N$-good. As a consequence, since 
${\rm d}(i, E \backslash F_M ) \geq M \geq N/2 $, the site $ i $ is  $(A,N)$-good. 
\end{itemize}
This concludes the proof of the Lemma.
\end{pf}

\section{Separation properties of the bad sites}\label{sec:sepa}\setcounter{equation}{0}

We now verify the ``separation properties" of the % $(A,N) $-
bad sites required in the multiscale Proposition \ref{propinv}. 

Let $ A := A(\e,\l,\teta ) $ be  the infinite dimensional matrix  of (\ref{matrpar}). 
We define % for $ j_0 \in  \Z^d \setminus {\cal Q}_N $, 
\be\label{tetabad}
B_{M}(j_0; \l )  :=  B_{M}(j_0; \e, \l ) :=  
\Big\{ \teta \in \R \, : \,   A_{M, j_0}(\e, \l,\theta) \ {\rm is} \ M-{\rm bad}  \Big\} \, . 
\ee
\begin{definition} {\bf ($ N$-good/bad parameters)} \label{def:freqgood}
A parameter $ \l \in \Lambda $ is $ N $-good for $ A $ if
\be\label{BNcomponents}
 \forall  M \in \{ N, N - 2 L_0 \} \, , \quad 
\forall \, j_0 \in \Z^d \setminus {\cal Q}_M \, , \quad 
B_{M}(j_0; \l) \subset \bigcup_{q = 1, \ldots, N^{2d+\nu+3}} I_q 
\ee
where $ I_q $ are  intervals with measure $ | I_q| \leq N^{-\t} $. 
Otherwise, we say $ \l  $ is $ N$-bad.
We define 
\be\label{good}
{\cal G}_{N} := {\cal G}_{N}({u}) := \Big\{ \l  \in  \Lambda \, : \,   \l  \ \  {\rm is \ } \ N-{\rm good \ for \ } A  \Big\} \, . 
\ee
\end{definition}

In order to prove the separation properties of the $ N $-bad sites
we have to require that $ \om = \l \bar \om $ satisfies a Diophantine type  non-resonance condition.
%(stronger than for the NLS equation). 
We assume:

\begin{itemize}
\item $ {\bf (NR)} $ There exist $ \g  > 0 $ such that, 
for any {\it non zero} polynomial $ P(X) \in \Z [X_1, \ldots, X_\nu] $ of the form 
\be\label{NRom0}
P( X ) = n + \sum_{1\leq i\leq j \leq \nu}  p_{ij} X_i X_j \, , \quad n, p_{ij} \in \Z \, , 
\ee  
we have
\be\label{NRom}
| P(\om)  | \geq \frac{\g}{ 1 +  |p|^{\t_0}} \, . % \quad \langle p \rangle := \max \{1, |p| \}  \, .  
\ee
\end{itemize}
The non-resonance condition $ ({\bf NR}) $ is satisfied by  $ \om = \l \bar \om $ for most $ \l \in \Lambda $, see 
Lemma \ref{sepadiof}.

\begin{remark}\label{rem:bou}
In  \cite{B5}, Bourgain  requires the  non-resonance condition 
\eqref{NRom}  for all non zero polynomials $ P(X) \in \Z [X_1, \ldots, X_\nu] $ of degree $ {\rm deg} P \leq 10 d $.
\end{remark}

The main result of this section is the following proposition.  It will
enable to verify the assumption (H3) of Proposition \ref{propinv} 
for the submatrices $ A_{N',j_0}(\e, \l,\theta) $. 

\begin{proposition}\label{prop:separation}
{\bf (Separation properties of $ N $-bad sites)}
There exists $ C_1(d, \nu) \geq 2 $, $N_0(\nu,d, \g_0, \Theta) \in \N$  such that
$ \forall N \geq N_0(\nu,d, \g_0, \Theta) $, if
\begin{itemize}
\item {\bf (i)}  $ \l $ is $ N $-good for $ A $,  % for all $ M \in \{ N - 2L_0, N \} $, 
\item {\bf (ii)} 
$ \tau > \chi \nu $,  % (= the  exponent of $ \bar \om $ in (\ref{diophan0})),
\item {\bf (iii)} 
$ \om = \l \bar \om $ satisfies $ {\bf (NR)} $, 
\end{itemize}
then, $  \forall \teta \in \R $,   
the $ N $-bad sites $ i:= (l,j) \in \Z^\nu \times \Z^d $  of $ A(\e, \l, \teta) $ with $|l| \leq N':=N^\chi $
admit a partition $ \cup_\a \Omega_\a $ in disjoint clusters satisfying
\be\label{separ}
{\rm diam}(\Omega_\a) 
\leq N^{C_1(d,\nu)} \, ,   \quad {\rm d}(\Omega_\a, \Omega_\b) > N^2 \, , \ 
\forall \a \neq	 \b \, .
\ee
\end{proposition}

The rest of this section is devoted to the proof of Proposition \ref{prop:separation}.
Note that, by (\ref{diophan0}), the frequency vectors $ \om = \l \bar \om $, $ \forall \l \in [1/2, 3/2] $, are Diophantine, namely
\be\label{diophan}
|\om \cdot l | \geq \frac{\g_0}{|l|^{\nu}} \, , \quad \forall l \in \Z^{\nu} \setminus \{ 0 \} \, . 
\ee 
\begin{lemma}\label{Ntime} Assume  that
$ \lambda $ is $ N $-good for $ A $  % for all $ M \in \{ N - 2L_0, N \} $
 %(i)-(ii) of Proposition \ref{prop:separation}. 
and let $\tau > \chi  \nu $. Then,  for all $ M \in \{ N - 2L_0, N \} $, 
$ \forall \bar{\jmath} \in \Z^d \backslash {\cal Q}_ M $, 
the number of $ M $-singular  sites $ (l_1,\bar{\jmath}) \in \Z^\nu \times \Z^d  $ 
with $ |l_1| \leq 2 N'  $  does not exceed $  N^{2d+\nu+3} $.
\end{lemma}

\begin{pf}
If $(l_1, \bar{\jmath}) $ is $ M $-singular then 
$ A_{M ,l_1,\bar{\jmath}}(\e, \l, \teta) $ is $ M $-bad (see Definitions \ref{GBsite} and \ref{goodmatrix} 
with $ N = M $). 
By the co-variance property (\ref{shifted}), we get that $ A_{M,\bar{\jmath}}(\e, \l, \teta + \l \bar{\om} \cdot l_1) $ is $ M $-bad, 
% and  $ |A_{(0,j_1)}^{(0,j_1)} (\e,\l,\theta + \l \bar{\om} \cdot l_1 )| \leq \Theta $ (or \eqref{tetaj0}), 
namely
$ \teta + \l \bar{\om} \cdot l_1 \in B_{M}(\bar{\jmath}; \l ) $, see (\ref{tetabad}). % with $ M = N $).
By assumption, $ \l$ is $ N$-good, and, therefore, (\ref{BNcomponents}) holds for $ M = N $ and $ M = N - 2 L_0 $.

We claim that in each interval $ I_q $ there is at most one element $ \theta + \om \cdot l_1  $ with 
$ \o = \l \bar{\o} $, $ |l_1| \leq 2N' $. 
Then, since there are at most $ N^{2d+\nu+3} $  intervals $ I_q $ (see (\ref{BNcomponents})),
  the lemma follows. 

We prove the previous claim by contradiction. 
Suppose that there exist $ l_1 \neq l_1' $ with $ |l_1|, |l_1'| \leq N'$, 
such that $ \om \cdot l_1 + \theta $, $ \om \cdot l_1' + \theta \in I_q $.  Then 
\be\label{up}
|\om \cdot (l_1 - l_1')| = |(\om \cdot l_1 + \theta) - (\om \cdot l_1' + \theta) | \leq |I_q| \leq N^{- \tau} \, .
\ee
By (\ref{diophan}) we also have
\be\label{lo}
|\om \cdot (l_1  - l_1')| \geq \frac{\g_0}{|l_1-l_1'|^{\nu}}  \geq \frac{\g_0}{ (4 N' )^{\nu}} = 4^{-\nu} \g_0 N^{- \chi \nu} \, . 
\ee
By assumption (ii) of Proposition \ref{prop:separation} 
the inequalities (\ref{up}) and (\ref{lo}) are in contradiction, for $ N \geq N_0(\g_0) $ large enough. 
\end{pf}

% We are now able to deduce an estimate of the $ N$-bad sites  $ (l_1,  \tilde \jmath  ) $ for all $\forall  \tilde \jmath  \in \Z^d $.

\begin{corollary}\label{Nbadtime}
Assume (i)-(ii)-(iii) of Proposition \ref{prop:separation}. Then, 
$ \forall  \tilde \jmath  \in \Z^d $, the number of $ N $-bad sites $ (l_1,  \tilde \jmath  ) \in \Z^\nu \times \Z^d  $ 
with $ |l_1| \leq N'  $  does not exceed $  N^{3d+2\nu+4} $. 
\end{corollary}

\begin{pf}
By Lemma \ref{Ntime}, for $M \in \{ N- 2L_0, N \} $, 
the set $ S_M $ of $M$-singular sites $ (l,j) \notin \check {\cal Q}_M $ (see \eqref{calQN} with $ N = M $)  with
$ |l| \leq N'+N$, $|j-   \tilde \jmath  | \leq M $ has cardinality at most
 $ C N^{2d+\nu +3} \times N^d $. 
Each $ N $-bad site $ (l_1,  \tilde \jmath  ) $ % $projecting to $j_1$
 with $|l_1| \leq N' $  is included, for some $M \in \{ N- 2L_0, N \} $,  in some
$ M $-ball  centered  at an element $ (l,j) $ of $ S_M $ which is not in $ \check {\cal Q}_M $ (see 
Definition \ref{GBsite}). 
Each of these balls contains at most $CN^\nu$ sites of the form $(l,  \tilde \jmath  ) $. Hence there are at most
$ C 2 N^{2d+\nu +3} \times N^d \times N^{\nu}$  such $N$-bad sites.
\end{pf}

We underline that the bound on the $ N $-bad sites given in  Corollary \ref{Nbadtime}
holds {\sl for all} $   \tilde \jmath  \in \Z^d $,
even if the complexity bound \eqref{BNcomponents} holds for all $ j_0 \notin {\cal Q}_M $. 
We now estimate  also the spatial components of the singular sites. 

\begin{definition} \label{chain}
{\bf ($ \Gamma $-chain)} 
A sequence $ i_0, \ldots , i_L \in \Z^{d+\nu}  $ of distinct integer vectors
satisfying 
$$ 
| i_{q+1} - i_q | \leq \Gamma \, , \quad \forall q = 0, \ldots, L - 1 \, ,
$$ 
 for some $ \Gamma  \geq 2 $,  is called a $ \Gamma $-chain of length $ L $.
\end{definition}

The next Lemma improves Lemma 20.14 of Bourgain \cite{B5}.
%requiring the weaker non-resonance assumption $ {\bf (NR)} $ (and a simpler proof). 

\begin{lemma}\label{Bourgain}
Assume that $ \om = \l \bar \om $ satisfies $  {\bf (NR)} $. 
For all $ \theta \in \R $, consider a $ \Gamma $-chain  $(l_q, j_q)_{q=0, \ldots, L} $ of $ \teta $-singular sites
with $  \Gamma \geq 2 $, namely,
$ \forall q = 0 , \ldots, L, $
\be\label{singusites}
\Big| (\l \bar \om \cdot l_q + \theta)^2 - \| j_q \|^2 - m  \Big| <  \Theta +1 \, , % \quad \Theta' := \Theta + m + 1 \, , 
\ee
such that, $ \forall \tilde{\jmath} \in \Z^d $, the cardinality
\be\label{cardinality}
|\{ (l_q,j_q)_{q = 0, \ldots, L} \, : \,  j_q = \tilde{\jmath} \}| \leq K \, .
\ee
Then its length is bounded by
\be\label{lunghez}
L \leq ( \Gamma K)^{C_2(d,\nu)} \, . 
\ee
\end{lemma}

\begin{pf} 
First note that  it  is sufficient to 
bound the length of a $ \Gamma $-chain of singular sites when $ \theta = 0 $. 
Indeed, suppose first that $ \theta = \o \cdot \bar l $ for some $ \bar l \in \Z^\nu $.
For a $ \Gamma $-chain of $ \theta $-singular sites $(l_q, j_q)_{q=0, \ldots, L} $, see \eqref{singusites},
the translated $ \Gamma $-chain
$(l_q + \bar l, j_q)_{q=0, \ldots, L} $, 
is formed by $ 0 $-singular sites, namely
$$
|(\om \cdot (l_q + \bar l))^2 - \| j_q \|^2 - m | < \Theta  \, . 
$$
For any $ \teta \in \R $, we consider an approximating sequence % it via a sequence of 
$ \o \cdot {\bar l}_n \to \teta $, $ {\bar l}_n \in \Z^\nu  $. 
A $ \Gamma $-chain of $ \theta $-singular sites (see \eqref{singusites}), is, for 
$ n $ large enough, also a $ \Gamma $-chain of $ \om \cdot {\bar l}_n $-sites. Then
 we bound its length arguing as in the above case.
\\[1mm]
We now introduce the quadratic  form 
$ Q : \R \times \R^d \to \R $ defined by
\be\label{Qxy}
Q(x,y) := - x^2 + \| y \|^2  
\ee
and the  associated bilinear symmetric form $ \varPhi : (\R \times \R^d)^2 \to \R $ defined by 
\be\label{def:varphi}
\varPhi \Big((x,y), (x',y') \Big) := - x x' + y \cdot y' \, .  
\ee
Note that $ \varPhi $ is the sum of the bilinear forms
\be\label{-v12}
\varPhi = - \varPhi_1 + \varPhi_2 
\ee
\be\label{vphi12}
\varPhi_1 \Big((x,y), (x',y') \Big) :=  x x' \, , \quad \varPhi_2 \Big((x,y), (x',y') \Big) := y \cdot y' \, . 
\ee
Let  $ (l_q, j_q)_{q= 0, \ldots, L} $ be a $ \Gamma $-chain,   namely 
\be\label{chainMlj}
| l_{q+1}- l_q|, | j_{q+1} - j_q | \leq \Gamma \, , \ \ \forall q=0, \ldots, L - 1 \, , 
\ee
of $ 0$-singular sites, see \eqref{singusites} with $ \teta = 0 $.
Setting 
\be \label{defxq}
x_q := \o \cdot l_q \, \in \, \o \cdot \Z^\nu \, ,
\ee
we get that (see \eqref{Qxy}) 
\be\label{singsi}
|Q ( x_q, j_{q})| < \Theta + 1+|m|  \, , \  \ \forall q = 0, \ldots, L  \, .
\ee
\begin{lemma}\label{T-+}  
$ \forall q, q_0 \in [0,L] $ we have
\be\label{seconda}
\Big| \varPhi \Big(( x_{q_0} , j_{q_0} ), ( x_q - x_{q_0}, j_q - j_{q_0}) \Big) \Big| \leq C |q - q_0 |^2 \Gamma^2 \, . 
\ee
\end{lemma}

\begin{pf}
By bilinearity 
\be\label{espbili}
Q ( x_q, j_q ) = Q ( x_{q_0}, j_{q_0} ) + 2 \varPhi \Big(( x_{q_0}, j_{q_0} ), ( x_q - x_{q_0}, j_q - j_{q_0}) \Big)  + 
Q ( x_q - x_{q_0},  j_q - j_{q_0} ) \, . 
\ee
We have
\begin{eqnarray}\label{boundch}
|Q ( x_q - x_{q_0},  j_q - j_{q_0} )| & \stackrel{\eqref{Qxy}} \leq & |x_q - x_{q_0} |^2 + \| j_q - j_{q_0} \|^2 \nonumber   \\
& \stackrel{\eqref{defxq}, \eqref{supeuc}}  
  \leq & |\om |^2 | l_q - l_{q_0}|^2 + d |j_q - j_{q_0}|^2 \stackrel{\eqref{chainMlj}}  \leq C |q-q_0|^2 \Gamma^2 \, . 
\end{eqnarray}
Then \eqref{seconda} follows by \eqref{espbili}, \eqref{boundch} and  \eqref{singsi}. 
\end{pf}

We introduce the subspace of $ \R^{d+1} $
\be\label{G+}
G := {\rm Span}_\R \Big\{  ( x_q - x_{q'}, j_q - j_{q'}) \, : \, 0 \leq q, q' \leq L  \Big\} = 
{\rm Span}_\R \Big\{ ( x_q - x_{q_0}, j_q - j_{q_0})  \, : \, 0 \leq q \leq L  \Big\}
\ee
and we call $ g \leq d + 1 $ the dimension of $ G $. Introducing a small parameter
$\delta >0$, to be specified later,we distinguish two cases. \\[1mm]
{\bf Case I}. 
$ \forall q_0 \in [0, L] $, 
\be\label{Case1} 
{\rm Span}_\R \{ ( x_q - x_{q_0}, j_q - j_{q_0}) \, : \,  | q - q_0 | \leq L^\d \, ,  \ q \in [0,L] \, \} = G \, . 
\ee
We select a basis of  $ G  \subset \R^{d+1}  $  
from  $ ( x_q - x_{q_0}, j_q - j_{q_0})  $ with $ | q - q_0 | \leq L^\d $, say 
\be\label{fs+}
f_s := ( x_{q_s} - x_{q_0}, j_{q_s} - j_{q_0}) = ( \om \cdot \Delta_sl , \Delta_s j ) \, , \ \ s = 1, \ldots, g \, ,
\ee
where
\be\label{lanuova}
(\Delta_s l, \Delta_s j) := ( l_{q_s} - l_{q_0}, j_{q_s} - j_{q_0} ) \quad {\rm satisfies} 
\quad  |(\Delta_s l, \Delta_s j)| \stackrel{\eqref{chainMlj}} \leq C \Gamma | q_s - q_0 | \leq C \Gamma L^\d \, . 
\ee
Hence
\be\label{boundfi}
| f_s | \leq C \, \Gamma L^\d \,  , \qquad \forall s =1, \ldots, g \, .
\ee

\begin{lemma}\label{Omegainv}
 Assume $ {\bf (NR)} $. 
Then the matrix 
\be\label{Omega}
\Omega := ( \Omega_s^{s'} )_{s,s'=1}^g \, , \quad \Omega_s^{s'} :=  \varPhi ( f_{s'}, f_s)  \, ,
\ee
is invertible and
\be\label{Cramer}
|( \Omega^{-1})_s^{s'}| \leq C (\Gamma L^\d)^{C_3(d,\nu)} \, , \quad \forall s, s' = 1, \ldots , g \, .
\ee
\end{lemma}

\begin{pf}
According to the splitting \eqref{-v12} we write  $ \Omega $  like
\be\label{-1-2}
\Omega :=  \Big(  - \varPhi_1 ( f_{s'}, f_s) + \varPhi_2 ( f_{s'}, f_s) \Big)_{s,s'=1, \ldots, g} = - S + R 
\ee
where, by \eqref{fs+}, 
\be\label{RSss'}
S_{s}^{s'} := \varPhi_1 ( f_{s'}, f_s) =
 (\o \cdot \Delta_{s'} l ) (\o \cdot \Delta_s l)   \, ,  
\quad R_{s}^{s'} := \varPhi_2 (f_{s'}, f_s) =  \Delta_{s'} j \cdot \Delta_s j  \, .
\ee
The matrix $ R = (R_1, \ldots, R_{g}) $ has integer entries (the $ R_i \in \Z^{g} $ denote the columns).
The matrix $ S := ( S_1, \ldots, S_g) $ has rank $ 1 $ since all its columns $  S_s \in \R^{g}  $ 
are colinear:
$$
S_s = (\om \cdot  \Delta_s l) ( \om \cdot \Delta_1 l , \ldots,   \om \cdot \Delta_g )^t \, , \quad s = 1, \ldots g \, . 
$$
We develop the determinant 
\begin{eqnarray}
P(\om) & := & {\rm det} \, \Omega  \stackrel{\eqref{-1-2}}  = {\rm det} (- S + R) \nonumber \\
&  = & {\rm det} (R ) - {\rm det} ( S_1, R_2, \ldots, R_{g}) - \ldots -
{\rm det} (R_1, \ldots, R_{g-1}, S_g) \label{polin}  
\end{eqnarray}
using that the determinant of matrices with $ 2 $ columns  $ S_{i} $, $ S_j $, $ i \neq j $,
is zero. The expression in \eqref{polin} is a polynomial in $ \om $ of degree $ 2 $ of the form \eqref{NRom0}
with coefficients 
\be\label{coeffb}
|(n,p)|  \stackrel{\eqref{RSss'}, \eqref{lanuova}} 
\leq C (\Gamma L^\d)^{C(d)} \, .
\ee 
If $ P \neq  0 $ then the non-resonance condition $ {\bf (NR)} $ implies % \eqrefm} and \eqref{coeffb} imply the lower bound 
\be\label{determinant}
|{\rm det} \, \Omega | = |P(\om)| \stackrel{ \eqref{NRom} } \geq \frac{\g}{1+|p|^{\t_0}} 
\stackrel{\eqref{coeffb} } \geq \frac{\gamma}{ (\Gamma L^\d)^{C'(d,\nu)}}  
\ee
(recall that $ \t_0 := \nu (\nu+1) $).
% which is \eqref{lowb-}. 
In order to conclude the proof of the lemma, we have to show that 
$ P \neq 0 $. % is not identically naught. 
By contradiction, if $ P= 0 $ then (compare with \eqref{-1-2})
% $ P( \ii \omega)  \equiv 0 $. But we have
$$
0 = P( \ii \omega)  = {\rm det} \Big( \varPhi_1 ( f_{s'}, f_s) + \varPhi_2 ( f_{s'}, f_s)\Big)_{s,s'=1,\ldots g} 
=  {\rm det} ( f_{s'} \cdot f_s)_{s,s'=1,\ldots g}   > 0
$$
because $ f_s  $ is a basis of $ \R^g $. This contradiction proves that $ P $ is not the zero polynomial. 

By \eqref{determinant}, the Cramer rule, and \eqref{boundfi} we deduce \eqref{Cramer}.
\end{pf}

We introduce 
$$
G^{\bot \varPhi} := \Big\{  z \in \R^{d+1}  \ : \  \varPhi (z, f) = 0 \, , \  \forall f \in G  \Big\} \, . 
$$
%Since $ \varPhi $ is non degenerate (see \eqref{def:varphi}) 
%$ {\rm dim} \, G + {\rm dim}  \, G^{\bot \varPhi}  = d+1 $ and,  
Since $ \Omega $ is invertible (Lemma \ref{Omegainv}), 
$\varPhi_{|G}$ is nondegenarate, hence
%we have
%$ G \cap G^{\bot \varPhi} = \{ 0 \} $. As a consequence 
$$
\R^{d+1} = G \oplus G^{\bot \varPhi} 
$$
and we denote by $ P_G : \R^{d+1} \to G $ the corresponding projector onto $ G $.

\smallskip

We are going to estimate 
\be\label{PGf}
P_G ( x_{q_0}, j_{q_0} ) = \sum_{s' =1}^{g} a_{s'} f_{s'} \, . 
\ee
For all $ s =1, \ldots, g $, and since $ f_s \in G $,  we have
$$
\varPhi \Big( ( x_{q_0}, j_{q_0} ), f_s \Big) = 
\varPhi \Big( P_G ( x_{q_0}, j_{q_0} ), f_s \Big) 
\stackrel{\eqref{PGf}} = \varPhi \Big( \sum_{s' =1}^{g} a_{s'} f_{s'}, f_s \Big) =
 \sum_{s' =1}^{g} a_{s'} \varPhi ( f_{s'}, f_s) 
$$
that we write as the linear system  
\be\label{def:ab}
\Omega a  = b \, , \qquad 
a := \left(\begin{array}{c} a_{1} \\ \ldots  \\ a_{g} \end{array}\right) \, , \quad 
b := \left(\begin{array}{c}  \varPhi \Big( ( x_{q_0}, j_{q_0} ), f_1 \Big)
\\ \ldots  \\  \varPhi \Big( ( x_{q_0}, j_{q_0} ), f_g \Big) \end{array}\right)
\ee
and  $ \Omega $ is defined in \eqref{Omega}.
\begin{lemma}\label{lem:12}
For all $ q_0 \in [0,L ] $ we have 
\be\label{PGs}
| P_G ( x_{q_0}, j_{q_0} ) | \leq ( \Gamma L^\d )^{C_4(d,\nu)} \, .
\ee
\end{lemma}

\begin{pf}
By \eqref{def:ab}, \eqref{fs+}, \eqref{seconda} and \eqref{Case1},   we get
%\be\label{boundb} 
$ |b| \leq C (\Gamma L^\d)^2 $. 
Hence, using also \eqref{def:ab} and \eqref{Cramer}, we get  %by   \eqref{Cramer}, we get
$ |a|  = | \Omega^{-1} b | \leq C (\Gamma L^\d)^{C}  $. 
This, with  \eqref{PGf} and \eqref{boundfi},  implies  \eqref{PGs}.
\end{pf}

As a consequence of Lemma \ref{lem:12}, for all $ q_1, q_2 \in [0,L] $, 
$$
| ( x_{q_1}, j_{q_1} ) -  ( x_{q_2}, j_{q_2} ) | =
| P_G \Big( ( x_{q_1}, j_{q_1} ) -  ( x_{q_2}, j_{q_2} )\Big) | 
\leq ( \Gamma L^\d )^{C_5(d,\nu)} \, .
$$
Therefore, for all $ q_1, q_2 \in [0,L] $, 
$ | j_{q_1} - j_{q_2} | \leq  (\Gamma L^\d)^{C_5(d,\nu)} $, 
and so
$$
{\rm diam}  \{ j_{q} \ ; \ 0\leq q \leq L \}  \leq   (\Gamma L^\d)^{C_5(d,\nu)} \, .
$$
Since all the $ j_{q}$ are in  $ \Z^d $, 
their number (counted without multiplicity) does not exceed $ C (\Gamma L^\d)^{C_5(d,\nu) d} $.
Thus we have obtained the bound
$$
\sharp \{ j_{q} \ : \ 0\leq q \leq L \} 
\leq C (\Gamma L^\d)^{C_5(d,\nu) d}  \, . 
$$
By assumption \eqref{cardinality}, for each $q_0 \in [0,L]$, the number of $q\in [0,L]$ such that
$j_q=j_{q_0}$  is  at most $ K $, 
and so
$$
L \leq  (\Gamma L^\d)^{C_6(d,\nu)}  K  \, . 
$$
Choosing $ \d > 0 $ such that $ \d C_6(d,\nu)  < 1 \slash 2 $, we get 
$ L \leq (\Gamma^{C_6(d,\nu) } K )^{2} $,
proving \eqref{lunghez}. 
\\[1mm]
{\bf Case II.} There is $ q_0 \in [0,L] $ such that  
$$
\mu := {\rm dim} \, {\rm Span}_\R \{  ( x_q - x_{q_0},  j_q - j_{q_0} ) \, : \, 
| q - q_0 | \leq L^\d \, ,  \ q \in [0,L] \, \}  \leq g -1 \, ,
$$
namely all the vectors $ (x_q, j_q) $ stay in a affine subspace of dimension $ \mu \leq g - 1 $. 
Then we repeat on the sub-chain $ (l_q, j_q)   $, $ | q - q_0 | \leq L^\d $, the argument of case I, 
to obtain a bound for $ L^\d $ (and hence for $L$).

Applying  at most $ (d+1) $-times the above procedure, we obtain a bound for $L$ of the form 
$ L \leq  ( \Gamma K)^{C(d,\nu)} $. This concludes the proof of Lemma \ref{Bourgain}. 
\end{pf}

\smallskip

\noindent 
{\sc Proof of Proposition \ref{prop:separation} completed.} 
Set  $ \Gamma :=  N^2  $ in Definition \ref{chain} and introduce the following equivalence relation: 

\begin{definition}\label{equivalence} 
We say that $ x \equiv y $ if  there is a $ N^2 $-chain $ \{ i_q \}_{q = 0, \ldots, L} $ of $ N $-bad sites
connecting $ x $ to $ y $, namely $ i_0 = x $, $ i_L = y $.
\end{definition}

A $  N^2 $-chain  $(l_q, j_q)_{q=0, \ldots, L} $ of $ N $-bad sites of $ A(\e, \l, \teta) $
 is formed by
$ \teta $-singular sites, namely \eqref{singusites} holds if $\e$ is small enough,  see Definition \ref{GBsite}.
Moreover, by Corollary \ref{Nbadtime} (remark it holds 
for all  $ \tilde \jmath \in \Z^\nu $), the condition \eqref{cardinality} of Lemma \ref{Bourgain}  is satisfied with 
$  K := N^{3d+ 2\nu+4} $. 
Hence Lemma \ref{Bourgain} implies %that the length 
\be\label{lenghtM}
L  \stackrel{\eqref{lunghez}} 
\leq ( N^2 N^{3d+ 2 \nu+4})^{C_2(d,\nu)} \leq N^{C'(d,\nu)} \, . 
\ee
The  equivalence relation in Definition \ref{equivalence} induces a partition of the $N$-bad sites 
of $ A(\e,\l, \teta) $ with $ |l| \leq N' $, 
in disjoint equivalent  classes $ ( \Omega_\a )$, satisfying  
$$
{\rm d}(\Omega_\a, \Omega_\b) >  N^2 \, , \quad {\rm diam}(\Omega_\a) \leq  N^2 L \stackrel{(\ref{lenghtM})} 
\leq N^2 N^{C'(d,\nu)} \leq N^{C_1(d,\nu)} \, .
$$

\section{Measure and complexity estimates}\label{sec:measure}\setcounter{equation}{0}

We define
\begin{eqnarray}\label{tetabadweak}
B_{N}^0 (j_0; \l) & := & B_{N}^0 (j_0; \e, \l) := \Big\{ \teta \in \R \, : \,   
\| A_{N, j_0}^{-1} (\e, \l,\theta)\|_0 > N^{\t}  \Big\} \\
& = &
\Big\{ \teta \in \R \, : \,  \exists {\rm \ an \ eigenvalue \ of \ }  A_{N, j_0} (\e, \l,\theta) \ {\rm with \ modulus \ less \ than} \ N^{-\t}  \Big\}
\label{la2}
\end{eqnarray}
where $ \|  \ \|_0 $ is the operatorial $ L^2 $-norm defined in (\ref{L2norm}).
The equivalence between (\ref{tetabadweak}) and (\ref{la2}) is a consequence of the 
self-adjointness of $A_{N, j_0} (\e, \l,\theta)$.  We also define
\begin{eqnarray}\label{weakgood}
{\cal G}_{N}^0 := {\cal G}_{N}^0 ({u}) & := & \Big\{ \l \in \Lambda \, : \,   \forall M \in \{ N, N - 2 L_0 \} \, , \
\forall \, j_0 \in \Z^d\backslash {\cal Q}_M  \, , \ 
B_{M}^0(j_0; \l) \subset \bigcup_{q = 1, \ldots, N^{2d+\nu+3}} I_q   \nonumber \\
& & \ \  {\rm where} \ I_q \ {\rm are \  \ intervals \ with \ measure} \ | I_q| \leq N^{-\t}   \Big\} \label{BNcomponent2}
\end{eqnarray}
(the set $ {\cal Q}_N $ is defined in \eqref{calQN}). 
The aim of this section is to provide, for any large $ N $, a suitable  bound % \eqref{measBN0}
on the Lebesgue measure of the complementary set of ${\cal G}_{N}^0 $, see \eqref{measBN0}. 
This will be used to estimate the measures of the sets
${\cal G}_{N}^c$ (see (\ref{good})) thanks to Proposition \ref{propinv}.

\begin{proposition}\label{PNmeas}
There are constants $ c, C > 0 $, $ N_0 \in \N $,  depending on $ V,d,\nu  $, such that, for all $ N \geq N_0 $ and
\be\label{ipopicco}
\e_0  (\| T_1 \|_0 + \| \partial_\l T_1 \|_0) \leq c 
\ee 
($ T_1 $  is defined  in \eqref{Tmatrix}$)$, 
the set $ {\cal B}_{N}^0 := \Lambda \setminus {\cal G}_{N}^0  $ 
has measure 
\be\label{measBN0}
|{\cal B}_{N}^0 | \leq C \,  N^{-1} \, .
\ee
\end{proposition}

The sequel of this section is devoted to the proof of Proposition \ref{PNmeas}. It is derived from several lemmas based on  
basic properties of eigenvalues of self-adjoint matrices,
which are a consequence of their variational characterization.  
In the definitions below, when $A$ is not invertible, we set $ \| A^{-1} \|_0 := \infty $. 

\begin{lemma}\label{variatione}% {\bf (see \cite{BB10})}
Let $J$ be an interval of $\R$ and $ A(\xi) $ be a family of self-adjoint
square matrices in ${\cal M}^E_E $, $ C^1 $ in the real parameter $ \xi  \in J $, and such that 
 $ \partial_\xi A(\xi) \geq \b I $ for some $ \b >  0  $. Then, for any $ \alpha >0 $, the Lebesgue measure 
$$
\Big| \Big\{ \xi \in J \, : \, \|A^{-1}(\xi)\|_0 \geq \alpha^{-1} \Big\}\Big| \leq  2 |E| \alpha \beta^{-1} 
$$
where $ | E | $ denotes the cardinality of the set $ E $. \\[1mm]
More precisely there is a family $(I_q)_{1\leq q \leq |E|}$ of intervals such that
\be |I_q| \leq 2\a \b^{-1} \quad {\rm and} \quad  
\Big\{ \xi \in J \, : \, \|A^{-1}(\xi)\|_0 \geq \alpha^{-1} \Big\} \subset 
\bigcup_{1\leq q \leq |E|} I_q
\ee
\end{lemma}

\begin{pf}
List the eigenvalues of the self-adjoint matrices $A(\xi)$ 
as $C^1$ functions $(\xi \mapsto \mu_q(\xi))$, $1\leq q \leq |E|$.  We have 
$$
\Big\{ \xi \in J \, : \, \|A^{-1}(\xi)\|_0 \geq \alpha^{-1} \Big\} =
\bigcup_{1\leq q \leq |E|} \Big\{ \xi \in J \, : \, \mu_q(\xi) \in [-\a , \a] \Big\}.
$$ 
Now, since $\partial_\xi A (\xi) \geq \beta I$, we have $\partial_\xi \mu_q (\xi) \geq \beta >0$, which
implies that $I_q := \{ \xi \in J \, : \, \mu_q(\xi) \in [-\a , \a] \}$ is an interval, of length less than 
$2\a \b^{-1}$. 
\end{pf}

\begin{lemma}\label{Lips}
Let $ A $, $ A_1 $ be self adjoint matrices. Then their eigenvalues (ranked in nondecreasing order)
satisfy the Lipschitz property
\be\label{continua}
|\mu_k (A) - \mu_k (A_1)| \leq \| A - A_1 \|_0  \, .
\ee
\end{lemma}

We develop all the computations for $ M = N $, the case $ M = N - 2 L_0 $ is the same.
We shall argue differently for $|j_0| \geq 8$ and $|j_0| < 8$ to estimate the complexity of
$B_N^0(j_0,\l)$. 

In the next lemmas we assume 
\be\label{assum1}
N \geq N_0 (V,\nu,d) > 0 \ {\rm large \ enough} \quad {\rm and} \quad 
 \e \| T_1 \|_0 \leq 1 \, .
\ee

%%%%%%%%%%%%%%%%%%%%%%%%%%%%%%%%%%%%%%%%%%%%%%%%%%%%%%%%

%%%%%%%%%%%%%%%%%%%%%%%%%%%%%%%%%%%%%%%%%%%%%%%%%%%%%%%%%%%%%%%%%%%%%%%%%%%%%%%%%%

\begin{lemma} \label{lem:complexity} 
$ \forall |j_0| \geq 8N  $, $ \forall \l \in \Lambda $,  we have
\be\label{complA}
B_{N}^0 (j_0;  \l) \subset \bigcup_{q=1, \ldots, 2(2N+1)^{d+\nu}} I_q 
\ee
where $ I_q $ are intervals  satisfying $  |I_q | \leq N^{- \t} $. 
%\label{sectio1} 
%$ \forall |j_0 | \geq 8 N $, $ \forall \l \in \Lambda $,  we have
%\be\label{B2alta}
%|B_{2,N}^0(j_0; \l)| \leq CN^{- \t +d +\nu-1} \, .
%\ee
\end{lemma}

\begin{pf}
We first claim that, if $ |j_0| \geq 8 N  $ and $  N \geq N_0(V,d,\nu)$ (see \eqref{assum1}), then  
\be\label{inclusionB}
B_{N}^0(j_0; \l) \subset \R \setminus [-4N,4N] \, . 
\ee
Indeed, by Lemma \ref{Lips} the eigenvalues  $ \lambda_{l,j}(\teta) $ of $ A_{N, j_0}(\e, \l,\teta)$ satisfy
\be \label{vpA}
\lambda_{l,j}(\teta) = \d_{l,j}(\teta)  + O(\e\| T_1 \|_{0} + \| V \|_0)  \quad {\rm where} \quad
\d_{l,j} (\teta) :=  - (\om \cdot l + \teta)^2 +  \|j\|^2  \, .
\ee
Since $ |\om | = |\l| | \bar \o | \leq 3/2 $ (see (\ref{baromega})), $\|j\| \geq |j|$ (see \eqref{supeuc}),  
$|j-j_0| \leq N $, $ |l| \leq N $, we get
\be\label{lo1}
\d_{l,j}(\teta)    \geq (|j_0| - |j-j_0|)^2   - (|\om||l| + |\theta| )^2 \geq 
(|j_0| - N)^2   -  (2N + |\teta|)^2 %\geq 10 N^2 
\ee
%since $ |j_0| \geq 8 N  $ and for  $ |\teta| \leq 4 N $.  
As a consequence, 
all the eigenvalues  $ \lambda_{l,j}(\teta) $  of $ A_{N, j_0}(\e, \l,\teta)$ satisfy, for 
$ |j_0| \geq 8 N  $ and  $ |\teta| \leq 4 N $, 
$$
\lambda_{l,j}(\teta) \stackrel{\eqref{vpA}, \eqref{lo1}} \geq 10 N^2 - O(\e\| T_1 \|_{0} + \|V \|_{0}) 
\stackrel{\eqref{assum1}} \geq N^2 \, , 
$$
implying \eqref{inclusionB}. We now estimate the complexity of 
$$
B_{N}^{0,-} :=  B_{N}^0(j_0;  \l) \cap (-\infty ,  - 4 N ) \quad {\rm and} \quad
B_{N}^{0,+} :=  B_{N}^0(j_0;  \l) \cap ( 4 N , \infty) \, .  
$$
Let us consider $ B_{N}^{0,-} $. 
For  $ \teta < - 4 N $,  the derivative 
$$
\partial_\teta A_{N, j_0}(\e, \l, \teta )  =
\, {\rm diag}_{|l| \leq N, |j - j_0| \leq N}  - 2( \om \cdot l + \teta) \\
>   8 N -  2|\om||l| \geq   5 N 
$$
and therefore Lemma \ref{variatione} (applied with $ \b = 5 N $, $ \a = N^{-\t} $) implies 
$$
B_{N}^{0,-} \cap (-\infty ,  - 4 N )  \subset  \bigcup_{1 \leq q \leq  (2N+1)^{d+\nu}} I^{-}_q , 
$$
where $ I_q^- $ are intervals  satisfying $  |I_q^- | \leq N^{- \t} $. 
We  get the same estimate for $  B_{N}^{0,+} $ and \eqref{complA} follows.
\end{pf}

We now consider the cases $ |j_0| < 8 N $.
Then the  continuity property (\ref{continua})  
of the eigenvalues allows to derive  a ``complexity estimate" for $ B_N^0 (j_0; \l) $
knowing its measure, more precisely the measure of  
\be\label{B2N0}
B_{2,N}^0(j_0;  \l) := B_{2,N}^0(j_0;  \e, \l) :=
\Big\{   \teta \in \R \, : \,  \| A_{N, j_0}^{-1} (\e, \l,\theta)\|_0 >  N^{\t} / 2  \Big\} \, .
\ee

\begin{lemma}\label{corol1}
$ \forall | j_0 | < 8 N $, $ \forall \l \in  \Lambda $, we have
\be\label{inclu0}
B_{2,N}^0(j_0;  \l) \subset I_N := [- 12 d N, 12  d N ] \, .
\ee
\end{lemma}

\begin{pf}
The eigenvalues $ \lambda_{l,j}(\teta)$ of $ A_{N,j_0}(\e,\l, \teta) $ satisfy \eqref{vpA} where, 
for  $ |\teta | \geq 12 d N $, 
\be\label{uppb}
|\om \cdot l + \teta| \geq |\teta | - |\om \cdot l  | \geq  12 d N - 2 N \geq 10 d N  \, ,
\ee
and, by \eqref{supeuc}, we have $ \|j\| ^2 \leq d (|j_0  | + |j-j_0| )^2 \leq d (9N)^2 $.  Hence
\begin{eqnarray*}
 \lambda_{l,j}(\teta) =  - (\om \cdot l + \teta)^2 +  \| j \|^2   + O(\e\| T_1 \|_0 + \|V \|_{0} ) 
&  \stackrel{\eqref{uppb}, \eqref{ipopicco}} \leq &
- (10dN)^2 + d (9N)^2   + C (1+\|V \|_0) \\
&  \leq & - 16 d^2 N^2
\end{eqnarray*}
for $ N \geq N(V,d,	\nu) $ large enough (see \eqref{assum1}), implying \eqref{inclu0}. 
\end{pf}

\begin{lemma}\label{complessita} %{\bf (see \cite{BB10})}
There is $ \hat{C} := \hat C (d) > 0 $  such that 
$ \forall |j_0| < 8 N  $, $ \forall \l \in  \Lambda $, we have 
$$ 
B^0_N (j_0;  \l) \subset \bigcup_{q=1,..., [\hat{C} \, {\mathtt M} N^{\tau +1} ] } I_q 
$$
where $ I_q $ are  intervals of length $ |I_q| \leq N^{-\tau } $ and  
$ {\mathtt M} := | B_{2,N}^0(j_0;  \l) |$.
\end{lemma}

\begin{pf}
Assume $\theta \in B_N^0 (j_0, \l)$, see \eqref{tetabadweak}. Then there is an eigenvalue of 
$A_{N,j_0}(\e, \l,\theta)$ with modulus less than $N^{-\t} $.  
Now, for $|\Delta \theta | \leq 1$, (recall \eqref{matrpar})
\begin{eqnarray*}
\| A_{N,j_0}(\e, \l,\theta + \Delta  \theta) - A_{N,j_0}(\e, \l,\theta) \|_0 & = &
\| {\rm Diag}_{|l| \leq N, |j-j_0| \leq N} \, (\l \ov{\om} \cdot l +\theta)^2 - (\l \ov{\om} \cdot l +\theta + \Delta \theta)^2 \|_0 \\
&\leq  &(4N+ 2| \theta | + 1 ) |\Delta \theta|.
\end{eqnarray*}
Hence, by Lemma \ref{Lips},   
\be\label{nuova}
(4N+ 2 | \theta | + 1 )  |\Delta \theta| \leq N^{-\t} \quad \Longrightarrow \quad  
\theta + \Delta \theta \in B_{2,N}^0 (j_0, \l)
\ee
because $ A_{N,j_0}(\e, \l,\theta + \Delta \theta)$
has an eigenvalue with modulus less than $2N^{-\t}$. %  {\it i.e.} $\theta + \Delta \theta \in B_{2,N}^0 (j_0, \l)$.
Now by Lemma \ref{corol1}, $|\theta | \leq 12dN$. Hence, by \eqref{nuova}, 
there is a positive constant $c := c(d) $ 
such that, for $\theta \in B^0_N (j_0;  \l)$, 
$$
%|\Delta \theta | \leq c N^{-(\t +1)} \Longrightarrow \theta + \Delta \theta \in 
[\theta - c N^{-(\t +1)} , \theta + c N^{-(\t +1)}] \subset  B_{2,N}^0 (j_0, \l) \, . 
$$ 
Therefore $ B_{N}^0(j_0, \l)  $ is included in an union of  intervals $J_m$ with disjoint interiors,
\be\label{inequa}
B_{N}^0 (j_0, \l)  \subset \bigcup_{m} J_m \subset B_{2,N}^0(j_0, \l), \quad {\rm with \ length} \quad  
|J_m|  \geq  2 c N^{-(\t+1)}
\ee
(if some of the intervals $ [\theta - cN^{-(\t+1)} , \theta +c N^{-(\t+1)} ]$ overlap, then we glue them together).  
We decompose each 
$ J_m $  as an union of (non overlapping) intervals $ I_q $ of length between 
$ c N^{-(\t+1)}/2 $ and  $ c N^{-(\tau +1)} $. Then, 
by (\ref{inequa}), we  get a new covering  
$$
B_{N}^0 (j_0, \l)  \subset \bigcup_{q=1, \ldots, Q} I_q  \subset B_{2,N}^0(j_0, \l) \quad  {\rm with} \ \   
c N^{-(\tau +1)} / 2 \leq |I_q | \leq c N^{-(\tau +1)} \leq N^{-\t}
$$
and, since the intervals $ I_q $ do not overlap, 
$$
Qc N^{- (\tau +1)}  \slash 2 \leq \sum_{q = 1}^Q | I_q |   
\leq |  B_{2,N}^0(j_0, \l) | =:   {\mathtt M} \, . 
$$
As a consequence $ Q \leq \hat C \,  {\mathtt M} \, N^{\t+1} $, proving the lemma. 
\end{pf}

The next lemma has major importance. 
%  the crucial change with respect to \cite{BB10}. 
The main difference with respect to the analogous lemma in
\cite{BB10} is that we do not assume the positivity of $ - \Delta + V(x) $, but only \eqref{nonke}.
% Then we perform a Lyapunov-Schmidt decomposition on the positive and negative directions.
Hence we have to require $ j_0 \notin {\cal Q}_N $.

\begin{lemma}\label{cor2} $ \forall |j_0| < 8 N $, $ j_0 \notin 
 {\cal Q}_N $, the set 
\be\label{B22N}
{\bf B}^{0}_{2,N}(j_0) := {\bf B}^{0}_{2,N}(j_0; \e) := 
\Big\{  (\l, \teta)  \in  \Lambda \times {\R} \, : \, 
\Big\| A_{N,  j_0}^{-1}(\e, \l,\teta) \Big\|_0 > N^{\t}/2  \Big\}
\ee
has measure 
\be\label{B2N} 
|{\bf B}^{0}_{2,N}(j_0)| \leq C N^{-\t+d+\nu+1}  \, . 
\ee 
\end{lemma}

\begin{pf}
By Lemma \ref{corol1},   $ {\bf B}^{0}_{2,N}(j_0)  \subset  \Lambda \times I_N $.
In order to estimate the  ``bad" $ (\l, \teta) $ where at least one eigenvalue of 
$ A_{N,j_0}(\e, \l, \teta )$ has modulus less than $ 2 N^{-\t} $, we introduce the variables
\be\label{changevaria}
\xi := \frac{1}{\l^2} \, , \quad \eta := \frac{\teta}{\l} 
\quad {\rm where} \quad (\xi, \eta) \in [4 \slash 9, 4] \times 2 I_N  \, .
\ee
Hence $ \teta = \l \eta $, $ \l := 1 / \sqrt{\xi} $, and we consider the self adjoint matrix
\be\label{etaxi}
A (\xi, \eta) := \frac{1}{\l^2} A_{N,j_0} (\e, \l,\teta) =
{\rm diag}_{|l| \leq N, |j - j_0| \leq N}\Big( - ( {\bar \om} \cdot l + \eta)^2 \Big)
+ \xi  P_{N,j_0} \, - \e \xi T_1(\e , 1 \slash \sqrt{\xi} )  
\ee
where, according to the notations \eqref{E0}-\eqref{DVB}, 
\be\label{PNpr}
P_{N, j_0} :=  \Pi_{N, j_0} ( - \Delta + V(x) )_{|E_{N,j_0}} \, .  % ( - \Delta + V(x) )_{| Q_N (j_0)}   
\ee
% with $ Q_N(j_0) := j_0 + [-N,N]^d $.  
% := \{ j \in \Z^d \, : \,  |j - j_0 | \leq N \}  $.
The  self-adjoint operator  $ P_{N, j_0} $ possesses a $ L^2 $-orthonormal basis of eigenvectors
$$
P_{N, j_0}  \Psi_j = \hat \mu_j \Psi_j 
$$
with real eigenvalues  $ (  \hat  \mu_j)_{j=1, \ldots (2N+1)^d}$ (depending on $ N $) 
indexed in non-decreasing order.
We define
$$
\quad  {\cal I}_- := \Big\{ j  
\ : \  \hat  \mu_j < 0   \Big\} \, , 
\qquad {\cal I}_+ := \Big\{ j \ : \  \hat  \mu_j > 0   \Big\} \, .
$$ 
Recalling the assumption $ j_0 \notin {\cal Q}_N $ (see \eqref{calQN}) 
% means that $ (0, \partial Q_N(j_0) ) \geq L_0 $ and
Lemma \ref{lzero} implies that: 
\begin{enumerate}
\item if $ {\mathtt B}(0,L_0-1) \subset \Z^d \setminus \{ |j - j_0| \leq N \} $ then $ P_{N, j_0} \geq \b_0 I $.
In this case  $ {\cal I}_- = \emptyset $,  $ {\cal I}_+ = \{ 1, \ldots, (2N+1)^d \} $ and  
$ \min_{j \in {\cal I}_+}  \hat  \mu_j \geq  \b_0 $.
\item  if $ {\mathtt B}(0,L_0) \subset \{ |j - j_0| \leq N \} $ then $ P_{N, j_0} $ has $ n^- $ % (independently of $ N $)
 negative eigenvalues $  \hat  \mu_j \leq - \b_0 $ and  the others $  \hat  \mu_j \geq \b_0 $
(we recall that $ n^- $ is the number of negative eigenvalues of $ -\Delta +V(x) $). We shall use that
 \be\label{massimo}
\max_{j \in {\cal I}_-}  \hat  \mu_j \leq - \b_0 \quad {\rm and} \quad \min_{j \in {\cal I}_+}  \hat  \mu_j \geq  \b_0  \, .
\ee
\end{enumerate}
% counted with multiplicity).  (in the above first case the set $ {\cal I}_- = \emptyset $) so that

We shall consider only the most difficult case 2 when % $ j_0 \notin {\cal Q}_N $,    $ 0 \in Q_N(j_0)  $ and so 
$ {\cal I}_- \neq \emptyset $. We  denote
$$
H_- := H_{{\cal I}_-} := \Big\{ u := \sum_{ |l | \leq N , j \in {\cal I}_-}  u_{l,j} 
e^{\ii l \cdot \vphi}  \Psi_j   \Big\} \, , \quad
H_+ :=  H_{{\cal I}_+} :=  \Big\{ u := \sum_{ |l| \leq N , j \in {\cal I}_+}  u_{l,j} 
e^{\ii l \cdot \vphi}  \Psi_j  \Big\} \, , 
$$
and  $ \Pi_- $, $ \Pi_+ $ the corresponding $ L^2 $-projectors.
Correspondingly we represent $ A := A(\xi, \eta) $ in \eqref{etaxi} as 
\be\label{Axi}
A = \left(\begin{array}{cc} 
A_- &  A_-^+ \\ 
A_+^- &  A_+
\end{array}\right) 
:=
 \left(\begin{array}{cc} 
 \Pi_- A_{| H_-}  &   \Pi_- A_{| H_+}  \\ 
 \Pi_+ A_{| H_-}  &   \Pi_+ A_{| H_+}
\end{array}\right) 
\ee
where $ A_-^+ = (A_+^-)^\dag $, $ A_-^\dag := A_- $ , $ A_+^\dag = A_+ $.

\begin{lemma}\label{A-1}
For all $ \xi \in [4/9, 4] $,  $ \eta \in \R $, % $ N \geq N_0 (V,d, \nu) $,
the matrix $ A_- := \Pi_- A_{| H_-} $ is invertible and
\be\label{A-1unif}
\| A_-^{-1} \|_0 \leq 3 \b_0^{-1} \, . 
\ee
\end{lemma}

\begin{pf}
By  \eqref{etaxi} and Lemma \ref{Lips}, 
the eigenvalues of the matrix $ A_-  $ satisfy, for $ |l| \leq N $, $ j \in {\cal I}_- $, 
\begin{eqnarray*}\label{eigneg}
 - ( \bar \om \cdot l + \eta)^2 + \xi  \hat \mu_j + O(\e \| T_1 \|_0) 
 & \leq &
 \xi \hat \mu_j + O(\e \| T_1 \|_0)  \leq \xi \max_{j \in {\cal I}_-} \hat \mu_j + O(\e  \| T_1 \|_0) \\
 & \stackrel{\eqref{massimo}, \eqref{ipopicco}} < &  -  \b_0 / 3 \, ,
\end{eqnarray*}
i.e. are negative and uniformly bounded away from zero. Then \eqref{A-1unif} follows. 
\end{pf}

The invertibility of the matrix in \eqref{Axi} is reduced to that of the self-adjoint matrix
\be\label{Lridu}
L := L(\xi, \eta):= A_+ - A_+^- A_-^{-1} A_-^+ 
\ee
via the  ``resolvent type" identity
\be\label{reso}
A^{-1} = \left(\begin{array}{cc} I & - A_-^{-1} A^+_- \\ 0 & I \end{array}\right)
\left(\begin{array}{cc} A_-^{-1} & 0 \\ 0 & L^{-1} \end{array}\right)
\left(\begin{array}{cc} I & 0 \\ - A^-_+ A_-^{-1} & I \end{array}\right) \, .
\ee
\begin{lemma}\label{L-1}
$\| L (\xi, \eta)^{-1} \|_0\leq N^\t / 20 $ except for  $ (\xi, \eta)  \in [4/9,4]  \times 2 I_N $ 
in a set of measure $ O( N^{- \t + d + \nu +1 })$. 
\end{lemma}

\begin{pf}
The derivative with respect to $ \xi $ of the matrix $ L(\xi, \eta) $ in (\ref{Lridu}) is
\begin{eqnarray}\label{lap1}
\partial_\xi L & = & \partial_\xi A_+ - (\partial_\xi A^-_+) A_-^{-1} A_-^+ - A^-_+ (\partial_\xi  A_-^{-1})  A_-^+ - 
A^-_+  A_-^{-1} ( \partial_\xi A_-^+ ) \nonumber \\
& = &  \partial_\xi A_+ - (\partial_\xi A^-_+) A_-^{-1} A_-^+  +  A^-_+  A_-^{-1} \, (\partial_\xi A_- ) A_-^{-1} A_-^+ - 
A^-_+  A_-^{-1} ( \partial_\xi A_-^+ ) \, .
\end{eqnarray}
Moreover, since $ \Pi_+ ( (\om \cdot \partial_\vphi)^2 - \Delta + V(x))_{|H_-} = 0 $ (and similarly exchanging $ \pm $), we have 
\be\label{Ala2} 
A_-^+ = - \e \xi \Pi_+ (T_1(\e , \xi^{-1/2}))_{| H_-} \, , \ \   A_+^- = - \e \xi \Pi_- (T_1 (\e , \xi^{-1/2}))_{| H_+} \, . 
\ee 
Hence, since $ 4 \geq \xi \geq 4/9 $, 
\be\label{precis}
\| A^+_- \|_0 + \| A^-_+ \|_0 + \| \partial_\xi A^+_- \|_0 + \| \partial_\xi A^-_+ \|_0 
= 0 ( \e (\| T_1 \|_0 + \| \partial_\l T_1 \|_0)).
\ee
In addition, by \eqref{etaxi}-\eqref{PNpr},
\be\label{ultim1}
\| \partial_\xi A_- \|_0 =  \| \Pi_- P_{N,j_0} |_{H_-} \|_0 +  O( \e ( \| T_1 \|_0 + \| \partial_\l T_1 \|_0) ) \leq C \, ,
\ee
\be\label{ultim2}
\partial_\xi A_+ = \Pi_+ P_{N,j_0}|_{H_+}  + 
O(  \e  (\| T_1 \|_0 + \| \partial_\l T_1 \|_0) \, .
\ee
Hence  by  \eqref{lap1},  \eqref{ultim2}, \eqref{precis}, \eqref{A-1unif},  \eqref{ultim1},  for $\e  (\| T_1 \|_0 + \| \partial_\l T_1 \|_0)$ small,
\be \label{varL}
\partial_\xi L      =   \Pi_+ P_{N,j_0}|_{H_+}  + 
O(  \e  (\| T_1 \|_0 + \| \partial_\l T_1 \|_0) 
\stackrel{\eqref{massimo}, \eqref{ipopicco}}  \geq  \frac{\b_0}{2} I \, .
\ee
 By (\ref{varL}) and Lemma \ref{variatione}, for each fixed 
$ \eta $, the set of $ \xi \in [4/9,4] $ such that at least one eigenvalue of the matrix $ L(\xi,\eta) $ in 
\eqref{Lridu} has modulus $ \leq  20 N^{-\t} $ has measure
at most $ O(  N^{- \t + d + \nu } \b_0^{-1})$. Then, integrating on $ \eta \in 2 I_N  $, whose length is $ |I_N| = O(N)$, 
we prove the lemma.
\end{pf}

From \eqref{reso},  \eqref{A-1unif}, \eqref{Ala2}, Lemma 
\ref{L-1} and
\eqref{ipopicco}, we derive the bound
\be\label{ilbou}
\| A^{-1} \|_0 \leq 
2 (\| L^{-1}(\xi, \eta) \|_0 + \| A_-^{-1} \|_0 ) \leq 2 \Big(  \frac{N^\t}{20} + 3 \b_0^{-1} \Big) 
\stackrel{\eqref{assum1}} \leq \frac{N^\t}{9}
\ee
except in a set of $ (\xi, \eta)  $ of measure $ O( N^{- \t + d + \nu +1 })$. 
We finally turn to the original parameters $ (\l,\teta) $.
Since the change of variables 
(\ref{changevaria})  has Jacobian of modulus greater than $ 1/8 $, 
% that the change of variables (\ref{changevaria}) has Jacobian  bounded  by $ 16 $. 
%see (\ref{B2N}), 
we have
$$ 
\| A_{N,j_0}^{-1}(\e,\l, \theta)\|_0 
\stackrel{\eqref{etaxi}}  = \l^{-2} \|A^{-1} \|_0 
\stackrel{ \eqref{baromega}, \eqref{ilbou}} \leq 4 \, \frac{N^\t}{9}  \leq  \frac{N^\t}{2} \, , 
$$
except for $ (\l, \theta)  \in \Lambda  \times \R $ in a set of measure $ \leq C N^{- \t + d + \nu + 1 } $.
The proof of Lemma \ref{cor2} is complete. 
\end{pf}

By the same arguments we also get the following measure estimate  used in the Nash-Moser iteration.

\begin{lemma}\label{measure0} 
The complementary of the set 
\be\label{Binver}
{\mathtt G}_{N} := {\mathtt G}_N({u}) := \Big\{ \l \in \Lambda \, : \, \| A_{N}^{-1}( \e, \l) \|_0 \leq  N^{\t}   \Big\}
\ee
%(the matrix $ A_N(\e,\l)$ is defined in \eqref{regulars})
has measure  
\be\label{measGN0} 
| \Lambda  \setminus {\mathtt G}_{N} | \leq  N^{- \t +d+\nu+2} \, . 
\ee 
\end{lemma}

As a consequence of Lemma \ref{cor2}, for ``most" $  \l$ the measure  of 
$ B_{2,N}^0 (j_0;  \l ) $ is ``small".

\begin{lemma}\label{intermed}
$ \forall |j_0| < 8 N $, $ j_0 \notin 
{\cal Q}_N $,  the set
$$
{\cal F}_{N}(j_0) :=  \Big\{ \l \in \Lambda \, : \, |B_{2,N}^0 (j_0;  \l)| \geq 
\hat{C}^{-1} N^{- \t + 2d + \nu + 2} \Big\} \, ,
$$
where $\hat{C}$ is the positive constant of Lemma \ref{complessita}, has measure 
\be\label{sectio}
| {\cal F}_{N}(j_0)| \leq C N^{- d -1 } \, .
\ee
\end{lemma}

\begin{pf}
By Fubini theorem (see \eqref{B22N} and \eqref{B2N0})
\be\label{Fubini}
|{\bf B}^{0}_{2,N}(j_0)| = \int_{\Lambda} | B_{2,N}^0(j_0;  \l) | \, d \l  \, .
\ee
Let  $ \mu := \t - 2 d - \nu - 2  $.  By (\ref{Fubini}) and (\ref{B2N}), 
\begin{eqnarray}
  C N^{-\t+d+ \nu + 1} & \geq & 
\int_{\Lambda} |B_{2,N}^0 (j_0;  \l) |  \, d \l \nonumber \\
& \geq & \hat{C}^{-1} N^{-\mu} \Big| \Big\{ \l \in \Lambda  \, : \, |B_{2,N}^0 (j_0; \l)| \geq 
\hat{C}^{-1} N^{-\mu} \Big\} \Big| \nonumber 
:= \hat{C}^{-1} N^{-\mu} |{\cal F}_{N}(j_0)| \nonumber
\end{eqnarray}
whence (\ref{sectio}).
\end{pf}

%By Lemma \ref{intermed}, 
For all $ \l \notin {\cal F}_{N}(j_0) $, % the measure 
$ |B_{2,N}^0(j_0; \l)| <  N^{- \t + 2d + \nu + 2} \hat{C}^{-1} $. 
Then Lemma \ref{complessita}  implies

\begin{corollary}\label{lem:complexity1} 
$ \forall |j_0| < 8N  $, $ j_0 \notin  {\cal Q}_N $, % \Z^d \setminus {\cal Q}_N $, 
$ \forall \l \notin {\cal F}_N(j_0) $, we have  
$$ 
B_{N}^0 (j_0;  \l) \subset \bigcup_{q=1, \ldots, N^{2d+\nu+3}} I_q 
$$ 
with $ I_q $ intervals  satisfying $ |I_q | \leq N^{- \t} $. 
\end{corollary}

Proposition \ref{PNmeas}  is now a direct consequence of the following lemma. 

\begin{lemma}\label{inclusion}
$ {\cal B}_{N}^0 \subseteq \bigcup_{|j_0| < 8 N, j_0 \notin  {\cal Q}_N} {\cal F}_{N}(j_0)  $.
\end{lemma}

\begin{pf}
Lemma \ref{lem:complexity}  and Corollary \ref{lem:complexity1} imply that
$$ 
\l \notin \bigcup_{|j_0| < 8 N,  j_0 \notin  {\cal Q}_N} {\cal F}_{N}(j_0)   \quad \Longrightarrow \quad 
\l \in  {\cal G}_N^0 
$$
(see the definition in (\ref{weakgood})).
The lemma follows. 
\end{pf}

\noindent
{\sc Proof of Proposition \ref{PNmeas} completed}.  
By Lemma \ref{inclusion} and (\ref{sectio}) we get
$$
| {\cal B}_{N}^0 | \leq \sum_{|j_0| < 8 N,  j_0 \notin  {\cal Q}_N} |{\cal F}_{N}(j_0)|  
\leq C (8N)^d    N^{-d-1} \leq  C  N^{-1} \, .
$$

\section{Nash Moser iterative scheme and proof of Theorem \ref{thm:main}}\label{sec:NM}\setcounter{equation}{0}

Consider the orthogonal splitting 
$$ 
H^s = H_n \oplus H_n^\bot
$$
where $  H^s $ is defined in \eqref{def:Hs} and 
$$ 
H_n :=  \Big\{ \, u = \sum_{|(l,j)| \leq N_n} \, u_{l,j}  \, e^{\ii (l \cdot \vphi + j \cdot x)}  
\Big\} \, ,  \quad 
H_n^\bot  :=  \Big\{  u = \sum_{|(l,j)|> N_n} \,   u_{l,j} \, e^{\ii (l \cdot \vphi + j \cdot x)} \in H^s \Big\} 
$$
with 
\be\label{defNn}
N_n := N_0^{2^n} \, ,  \ \qquad   {\rm namely} \ \ \ N_{n+1} = N_n^2  \, ,  \ \forall n \geq 0 \, . 
\ee
We shall take  $ N_0 \in \N $ large enough depending on $ \e_0 $ and $ V $, $ d $, $ \nu $.
Moreover we  always assume  $ N_0 > L_0 $ defined in Lemma \ref{lzero}.  
We denote by 
$$
P_n : H^s \to H_n \qquad \qquad {\rm and} \qquad \qquad   
P^{\bot}_n: H^s \to H_n^\bot
$$
the orthogonal projectors onto $ H_n $ and $ H_n^{\bot} $. 
The following ``smoothing" properties hold, $ \forall n \in \N $, $ s \geq 0 $, $ r \geq 0 $,  
$$
 \| P_n u \|_{s+r}  \leq  N_n^r \| u \|_s \, , \   \forall u \in H^s \, , \qquad % \label{S1}
 \|  P_n^\bot u \|_{s}  \leq   N_n^{- r } \| u \|_{s+r} \, , \   \forall u \in H^{s+r} \, . % \label{S2}
$$
For $ f \in C^q (\T^\nu \times \T^d \times \R; \R)$ with  
\be\label{defk} 
q \geq S + 2 \, ,
\ee
the  composition operator on  Sobolev spaces
$$
f:  H^s \to H^s \, , \qquad f(u) (\vphi,x) := f(\vphi,x,u (\vphi,x) )   \,  
$$
satisfies the following standard properties: % (see e.g.  \cite{LM})
 $ \forall s \in [s_1, S] $,
$ s_1 > (d+\nu)/2 $, 
\begin{itemize}
\item {\bf (F1)} ({\bf Regularity}) 
$ f \in C^2 ( H^{s}; H^{s} )$. 
\end{itemize}
\begin{itemize}
\item {\bf (F2)} ({\bf Tame estimates}) 
$ \forall u, h \in H^{s} $ with $ \| u \|_{s_1} \leq 1 $, 
\be\label{fDftame} 
\| f(u) \|_{s}  \leq C(s) (1 + \| u \|_s) \, , \ \  \| (Df)(u) h \|_s 
\leq C(s) ( \| h \|_s + \| u \|_{s} \| h \|_{s_1}) \, ,
\ee
$$
\| D^2 f(u)[h,v] \|_{s} \leq C(s) \Big(\|u\|_{s} \|h\|_{s_1} \|v\|_{s_1} + \|v\|_{s} \|h\|_{s_1} + \|v\|_{s_1} \|h\|_{s} \Big) 	\, . 
$$
\item {\bf (F3)} 
({\bf Taylor Tame estimate}) 
$ \forall u \in H^{s} $ with $ \| u \|_{s_1} \leq 1 $,  $ \forall h  \in H^{s} $ with $\| h \|_{s_1} \leq 1$, 
$$
\| f( u + h) - f (u) -  (Df)(u) \, h \|_{s}  \leq  C(s) ( \| u\|_{s}  \| h \|_{s_1}^2 + \| h \|_{s_1}   \| h \|_{s}) \, .
$$
In particular, for $ s = s_1 $, 
$ \| f (u+ h) - f (u) - (Df)(u) \, h \|_{s_1} \leq C(s_1) \| h \|_{s_1}^2 $.
\end{itemize}
We fix the Sobolev indices
$ s_0 < s_1 < S  $
as
\be\label{Sgr}
s_0 := b = d + \nu \, , \qquad  
s_1 := 10 (\t + b)C_2 \, , \qquad 
S := 12 \t' +  8(s_1 + 1)   \, ,
\ee
where  
\be\label{tautau0}
C_2 :=  6(C_1 + 2)  \,  , \ \t :=   \max\{ d + \nu + 3, 2  C_2   \nu + 1 \} \, ,
   \  \t' := 3 \t+ 2 b \, ,  % \tau_0 := \nu (\nu+1)  
\ee
% (the constant $ \tau_0   $ is introduced in \eqref{diofgr})
and $ C_1 := C_1(d,\nu) \geq 2 $ is defined in Proposition \ref{prop:separation}.
Note that $ s_0, s_1, S $ defined in \eqref{Sgr} 
depend only on $ d $ and $ \nu $. We  also fix the constant $ \d $  in Definition \ref{goodmatrix} as
\be\label{delta14}
\d := 1/ 4  \, .
\ee

\begin{remark}\label{defchi}
By (\ref{Sgr})-(\ref{delta14})  the hypotheses (\ref{dtC})-(\ref{s1}) of Proposition  \ref{propinv} 
are satisfied for any $ \chi \in [C_2, 2 C_2) $, 
as well as assumption (ii) of Proposition \ref{prop:separation}. 
We assume $ \t \geq d + \nu + 3  $ in view of  (\ref{measGN0}). 
% The strongest condition for $ S $ appears in the proof of  Lemma  \ref{deriva}.
\end{remark}

Setting 
\be\label{tau1}
\t_1 :=  3 \nu + d + 1 
\ee
and $ \gamma > 0 $, we  implement the first steps 
of the Nash-Moser iteration restricting $ \l $ to the set
\begin{eqnarray}\label{diofs}
\bar {\cal G} & := &  \Big\{ \l  \in \Lambda \, : \, 
\Big\| \Big( - \l^2 ( \bar \om \cdot l)^2 + \Pi_{0} (-\D+ V(x))_{| E_{0}}\Big)^{-1} \Big\|_{L^2_x} 
\leq \frac{N_0^{\t_1}}{\g}, 
 \,  \forall  \, |l| \leq N_0   \Big\} \nonumber \\
& = & 
\Big\{ \l  \in \Lambda \, : \, 
|  - \l^2 ( \bar \om \cdot l)^2 +  \hat \mu_j | \geq \g N_0^{-\t_1}, 
\, \forall  \, |j| \leq N_0, \, |l| \leq N_0    \Big\} 
\end{eqnarray}
where $ \hat \mu_j $ are the eigenvalues of $ \Pi_{0} (- \D + V(x) )_{| E_{0}} $ and
$ \Pi_0 := \Pi_{N_0 ,0} $, $ E_0 := E_{N_0 ,0} $ are defined in  \eqref{E0}.
We shall prove in Lemma \ref{NRfre}  that  $ |\bar {\cal G}| = 1 - O(\g) $ (since $ \t_1 >  3 \nu + d $).

%Setting 
%\be\label{tau0'}
%\t_0' := 2 + \nu(\nu+1) \, ,  %  \nu (\nu+1) + 1 \, , 
%\ee
We prove the separation properties of the small divisors for $ \l $ 
satisfying assumption (NR), namely  in % at each step on 
\begin{eqnarray}\label{ditilde}
\tilde {\cal G} & := & 
\Big\{ \l  \in \Lambda \, : \, 
\Big| n +  \l^2 \sum_{1 \leq i \leq j \leq \nu} p_{ij}  \bar \om_i \bar \om_j \Big| \geq \frac{\g}{1 + |p|^{\t_0}}
 \, , \ \  \forall  \, (n,p) \neq 0 \Big\} \, . 
\end{eqnarray}
The constant $ \g $ will be fixed in \eqref{N0ge}. We also set
\be\label{def:sigma}
\s := \t' + \d s_1 + 2   \, . 
\ee
Given a set $ A $ we denote $ {\cal N}(A, \eta) $ the open neighborhood of $ A $ of width $ \eta $ (which is empty if $ A $
is empty).

\begin{theorem}\label{cor1} {\bf (Nash-Moser)} 
There exist $\e_0 , \bar c , \bar \gamma > 0 $ {\rm (}depending on $ d,\nu,V$,$\g_0 )$ such that, if 
\be\label{smallsto}
\gamma \in (0, \bar \g)  \, , \  N_0 \geq 2 \g^{-1}  \, ,  
\qquad { and } \qquad  
\e \in [0, \e_0) \, , \ \e N_0^{S} \leq \bar c  \, , 
\ee
then there is a sequence $(u_n)_{n \geq 0} $ of $ C^1$ maps 
$ u_n (\e, \cdot) : $ $ \Lambda \to   H^{s_1} $ satisfying 
\begin{itemize}
\item[$ {\bf (S1)}_n $] \ $ u_n (\e, \l ) \in H_n  $,  $ u_n ( 0, \l ) = 0 $, 
$ \| u_n \|_{s_1} \leq 1$, $ \| u_0 \|_{s_1} \leq N_0^{-\s} $  and $ \| \partial_{\l} u_n \|_{s_1} \leq C(s_1) N_0^{\t_1+ s_1+ 1} \g^{-1} $.
\item[$ {\bf (S2)}_n $]  \ $ (n \geq 1) $ \  For all  $ 1 \leq k \leq n $, 
$\|u_k-u_{k-1}\|_{s_1} \leq N_k^{-\s -1 } $, 
$ \|\partial_{\l}(u_k-u_{k-1})\|_{s_1} \leq N_k^{-1/2}$. 
\item[$ {\bf (S3)}_n $] \  $ (n \geq 1) $ 
\be\label{GN0N}
\| u - u_{n-1} \|_{s_1} \leq N_n^{-\sigma} \quad  \Longrightarrow \quad
\bigcap_{k= 1}^n {\cal G}_{N_k}^0 (u_{k-1}) \cap   \tilde {\cal G}  
  \subseteq  {\cal G}_{N_{n}}(u) 
\ee
where  
$ {\cal G}^0_{N } (u ) $ (resp. ${\cal G}_{N } (u ) $) is defined in  (\ref{weakgood}) (resp. in (\ref{good}))  
and $ \tilde {\cal G}  $ in \eqref{ditilde}. 
\item[$ {\bf (S4)}_n $]
Define the set
\be\label{Gscavo}
{\cal C}_n := 
\bigcap_{k= 1}^n {\mathtt G}_{N_k}(u_{k-1})  \bigcap_{k= 1}^n 
{\cal G}_{N_k}^0(u_{k-1})   \bigcap {\tilde {\cal G}}  \cap  {\bar {\cal G}}  \,  , 
\ee
where  ${\mathtt G}_{N_{k}}(u_{k-1}) $ is defined in (\ref{Binver}), $ \bar {\cal G} $ in (\ref{diofs}),
$ \tilde {\cal G}  $ in \eqref{ditilde},
$ {\cal G}_{N_k}^0(u_{k-1}) $ in (\ref{weakgood}). 
 
If $ \l \in {\cal N} ( {\cal C}_n, N_n^{-\s}) $ then  
$  u_n(\e,\l) $ solves the equation
$$
P_{n} \Big(L_{\om} u - \e f( u )  \Big)= 0  \, . \leqno{(P_n)}   
$$
\item[$ {\bf (S5)}_n $] \   $ U_n :=  \|  u_n \|_{S} $,  
$ U_n' :=  \|  \partial_{\l}{ u}_n \|_S $  (where $ S $ is defined in (\ref{Sgr}))
satisfy
$$
(i) \;\; U_n \leq  \displaystyle  N_n^{ 2( \t' + \d s_1 + 1)} \, , 
\qquad  (ii) \;\; U_n' \leq  N_{n}^{4\t'+ 2 s_1 + 4}  \, . 
$$
\end{itemize}
The sequence $ (u_n)_{n \geq 0} $  converges 
in $ C^1 $ norm to a map
\be\label{uC1} 
u (\e, \cdot )  \in C^1 ( \Lambda ,  H^{s_1}) \quad {\rm with} \quad u(0, \l) = 0 
\ee
and,  if $ \l $ belongs to  the Cantor like set
\be\label{Cinfty}
{\cal C}_\e := 
\bigcap_{n \geq 0} {\cal C}_n 
\ee
then $ u(\e,\l) $  is a solution of  \eqref{eq:freq}, with $ \om = \l \bar \om $.
\end{theorem}

The proof of  Theorem \ref{cor1} follows exactly the steps in  \cite{BB10}, section 7.
A difference is that  we do not need to estimate  $ \partial_\e u_n $.  
Another difference is that 
the frequencies in $ {\cal C}_{n} $ (see \eqref{Gscavo}) belong also to $ \tilde {\cal G} $ (in order to prove the
separation properties). 
For the reader convenience, in the Appendix, we spell out the main steps indicating the other minor  adaptations in the proof.
The main one concerns the  proof of Lemma
\ref{S3n+1} where we 
estimate  $ A_{M, j_0}^{-1}(\e,\l,\teta ) $ for both $ M = N_{n+1} $ and $ N_{n+1} - 2 L_0 $ (and not only $ N_{n+1} $).

\smallskip 

The sets of parameters $ {\cal C}_n $  in $ (S4)_n $ are decreasing, i.e.
$$
\ldots \subseteq {\cal C}_{n} \subseteq {\cal C}_{n-1} \subseteq \ldots \subseteq  {\cal C}_{0} \subset 
 \tilde {\cal G} \cap  \bar {\cal G} \subset \Lambda \, , 
$$
and it could happen that $ {\cal C}_{n_0} = \emptyset $ for some $ n_0 \geq 1 $. 
In such a case $ u_n = u_{n_0} $, $ \forall n \geq n_0 $
(however the map $ u(\e, \cdot ) $ in \eqref{uC1} is always defined), and $ {\cal C}_\e = \emptyset  $.
%Later, in (\ref{N0ge}), we shall specify the values of $ \g, \e_0, N_0 $, in order to
We shall prove in \eqref{Cefull} that  (choosing \eqref{N0ge}) the set 
$ {\cal C}_\e $ has asymptotically full measure.

\smallskip

In order to prove Theorem \ref{thm:main}, 
we first verify the existence of frequencies satisfying \eqref{diofp}.

\begin{lemma}\label{diofp}
For  $ \t_0 > \nu (\nu+1) -1 $, the  complementary of the set of $ \om \in {\R^\nu} $, $ |\om| \leq 1 $,  verifying \eqref{diofgr} 
has measure $ O( \g_0^{1/2}) $.  
\end{lemma}

\begin{pf} 
We have to estimate the measure of 
$$
\bigcup_{p \in \Z^{\nu (\nu+1) /2} \setminus \{0 \}}  \!\! \!\! {\cal R}_p 
\quad {\rm where} \quad {\cal R}_p  := 
\Big\{ \omega \in \R^\nu \, , \, |\om | \leq 1  \ : \  
\Big| \sum_{1 \leq i \leq j \leq \nu} \omega_i \omega_j p_{ij} \Big| < \frac{\gamma_0}{|p|^{\tau_0}} \Big\} \, . 
$$
Let $ M := M_p $ be the $ (\nu \times \nu)-${\it symmetric} matrix such that 
$$
\sum_{1 \leq i \leq j \leq \nu} \omega_i \omega_j p_{ij}  = M \om \cdot \om \, , \quad \forall \om \in \R^\nu \, .
$$
The symmetric matrix $ M $ has coefficients   
\be\label{defMS}
M_{ij} := \frac{p_{ij}}{2} (1 + \d_{ij})  \, , \  \forall   1 \leq i \leq j \leq \nu \, , \quad {\rm and} \quad M_{ij} = M_{ji} \, .
\ee
There is an orthonormal basis of eigenvectors 
$ V := (v_1, \ldots, v_k) $ of  $ M v_k = \l_k v_k $ with real eigenvalues
$ \l_k :=  \lambda_k (p) $. 
Under the isometric change of variables $ \omega = V y $ we have to estimate 
\be\label{Rpy}
|{\cal R}_p | = \Big| \Big\{ y \in \R^\nu \, , \, | y | \leq 1  \ : \  
\Big| \ \sum_{1 \leq k  \leq \nu} \l_k y_k^2 \Big| < \frac{\gamma_0}{|p|^{\tau_0}} \Big\}\Big| \, . 
\ee
Since
$ M^2 v_k = \l_k^2 v_k $, $ \forall k =1, \ldots, \nu $,  we get 
$$
\sum_{k=1}^\nu \l_k^2 = {\rm Tr} (M^2) = \sum_{i, j =1}^\nu M_{ij}^2 
\stackrel{\eqref{defMS}} \geq | p |^2 / 2  \, . 
$$
Hence there is an index $ k_0 \in \{1, \ldots , \nu \} $ 
such that $ | \l_{k_0} | \geq  |p| / \sqrt{2\nu} $ and  the derivative
\be\label{Rpy1}
\Big| \partial^2_{y_{k_0} } \Big( \sum_{1 \leq i  \leq \nu} \l_k y_k^2 \Big) \Big| =   |2 \l_{k_0}| \geq 
\sqrt{2} \, |p| \slash \sqrt{\nu}  \, . 
\ee
As a consequence of \eqref{Rpy} and \eqref{Rpy1} we deduce the measure estimate
$ |{\cal R}_p | \leq  C  \sqrt{\frac{\gamma_0}{|p|^{\tau_0 +1}} } $ (see e.g. Lemma 9.1 in \cite{EK})
and 
$$
\Big| \bigcup_{p \in \Z^{\nu (\nu+1) /2} \setminus \{0 \}}  \!\! \!\! {\cal R}_p  \Big| \leq 
 \sum_{p \in \Z^{\nu (\nu+1) /2 } \setminus \{0\}}  \!\! \!\!  |{\cal R}_p | \leq  \sum_{p \in \Z^{\nu (\nu+1) /2}
 \setminus \{ 0 \}}  \!\! \!\!    C  \sqrt{\frac{\gamma_0}{|p|^{\tau_0 +1}} } \leq C' \sqrt{\gamma_0}
$$
for $ \tau_0 > \nu (\nu +1 ) - 1 $. 
\end{pf}

\smallskip

We now  prove that $ {\cal C}_\e  $ in \eqref{Cinfty} has asymptotically full measure, i.e. (\ref{Cmeas}) holds.

\begin{lemma}\label{NRfre} 
The complementary of the set  $ \bar {\cal G} $ defined in \eqref{diofs} satisfies
\be\label{calGg}
| \Lambda \setminus \bar {\cal G} | =  O(\g) \, . 
\ee
\end{lemma}

\begin{pf}
The  $ \l $ such that  (\ref{diofs}) is violated are 
\be\label{eq:1}
\Lambda \setminus {\bar {\cal G}} = % \subseteq 
\bigcup_{|l|,|j| \leq N_0} {\cal R}_{l,j} \quad 
{\rm where} \quad 
{\cal R}_{l,j} := \Big\{ \l \in [1/2, 3/2 ] \, : 
\, | \lambda^2 ( \bar \om \cdot l)^2 - \hat \mu_j | < \frac{\g}{N_0^{\t_1}}  \Big\} \, . 
\ee
By Lemma \ref{lzero} the eigenvalues $ | \hat \mu_j | > \b_0 $ (for $ N_0 > L_0 $ ). Therefore, 
$ {\cal R}_{0,j} = \emptyset $ if $ \g N_0^{-\t_1} < \b_0 $.
We have to estimate the $ \xi :=  \l^2  \in [4/9,4] $ such that 
$ |  \xi (\bar \om \cdot l)^2 -  \hat \mu_j | <  \g N_0^{-\t_1} $. 
The derivative of
the function $ g_{lj} (\xi) :=   \xi (\bar \om \cdot l)^2 -  \hat \mu_j  $ satisfies
$ \partial_\xi g_{lj} (\xi) =  (\bar \om \cdot l)^2  \geq 4 \g_0^2 N_0^{ - 2 \nu } $ by  \eqref{diophan0}.
%because $ \Pi_0 (- \D + V(x) )_{| E_0} \geq \b_0 I   by \eqref{eq:posi}. 
As a consequence
\be\label{eq:3}
|{\cal R}_{l,j}| \leq  C  \g_0^{-2} \g N_0^{-\tau_1 + 2 \nu } \, . 
\ee
Then (\ref{eq:1}), (\ref{eq:3}), imply
$$
| \Lambda \setminus {\bar {\cal G}} | \leq \sum_{|l| \leq N_0, |j| \leq  N_0} |{\cal R}_{l,j}| \leq
C  \g \g_0^{-2} N_0^{d+\nu} N_0^{- \tau_1 + 2 \nu }= O(\g)
$$
since $ \t_1  >  3 \nu + d  $ (see \eqref{tau1}).
\end{pf}

\begin{lemma}\label{sepadiof}
Let $ \g \in (0, 1/4) $. Then the complementary of the set 
$ \tilde {\cal G} $   
in \eqref{ditilde}  has a measure  
\be\label{tildeGM} 
|\Lambda \setminus \tilde {\cal G} | = O(\gamma  ) \, .
\ee
\end{lemma}

\begin{pf}
For $ p := (p_{ij})_{1 \leq i \leq j \leq \nu} \in \Z^{\nu(\nu +1)/2} $, let 
$$
a_p := \sum_{1\leq i\leq j \leq \nu} p_{ij} {\bar \om}_i {\bar \om}_j \,  , \quad  g_{n,p}(\xi) :=  n + \xi a_p \, . 
$$
We have 
\be\label{somma}
|\Lambda \setminus \tilde {\cal G} |  \leq 
C \sum_{(n,p) \neq (0,0)} |{\cal R}_{n,p}|  \quad  {\rm where} 
 \quad  {\cal R}_{n,p} := \Big\{  \xi := \l^2 \in  [1/4, 9/4] \, : \,   | g_{n,p}(\xi) | <    \frac{\g}{ 1 + |p |^{\t_0}} \Big\}
\ee
%and 
%$$
%{\cal R}_{n,p} := \Big\{  \xi := \l^2 \in  [1/4, 9/4] \, : \,   | g_{n,p}(\xi) | <    \frac{\g}{ 1 + |p |^{\t_0}} \, , 
%\  g_{n,p}(\xi) :=  n + \xi \sum_{1\leq i\leq j \leq \nu}  
%p_{ij} {\bar \om}_i {\bar \om}_j   \Big\} \, . 
%$$
%The sum \eqref{somma} can be limited to the indices $ |n | \leq C |p| $, otherwise $ {\cal R}_{n,p} = \emptyset $.
\\[1mm]
{\bf Case I:} $ n \neq 0 $.  If $ {\cal R}_{n,p} \neq \emptyset $ then, since $ \g  \in (0, 1/4) $ and $|\xi| \leq 3$, we deduce  
$|a_p| \geq 1/4$, $|n| \leq 4 |a_p|$ and
%$$
%\Big| \sum_{1\leq i\leq j \leq \nu}  p_{ij} {\bar \om}_i {\bar \om}_j \Big| \geq  \frac14 \, , 
%$$
%namely $ g_{n,p}(\xi) := $ $  |\partial_\xi g_{n,p}(\xi)| \geq 1 / 4 $, and so 
$$
|{\cal R}_{n,p}| \leq  \frac{2 \g}{ (1 + |p |^{\t_0}) |a_p|}  \, . 
$$
Hence 
\be\label{se1}
\sum_{n \in \Z \setminus \{0\}} |{\cal R}_{n,p}| = \!\!
\sum_{n \in \Z \setminus \{0\}, |n| \leq 4|a_p|} \!\! |{\cal R}_{n,p}| \leq \frac{C \g}{(1+|p|)^{\tau_0}}.
\ee
{\bf Case II:} $ n = 0 $. In this case, %each $ {\cal R}_{0,p} $, $ p \neq 0 $, 
%is an interval centered at the origin and,  
using \eqref{diofgr} we obtain 
\be\label{se2}
 {\cal R}_{0,p} \subset  \Big( 0, \frac{\gamma}{ 1 + |p |^{\t_0}} \frac{ |p|^{\t_0} }{\g_0} \Big] 
 \subset  \Big( 0, \frac{\gamma}{\g_0} \Big] \, .  % \quad \forall p \, . 
\ee
From \eqref{somma}, \eqref{se1}, \eqref{se2}, 
%  \eqref{tau0'},  
$ \t_0 := \nu (\nu+1) $, % > 1 + (\nu(\nu+1)/2) $ 
we deduce \eqref{tildeGM}. 
\end{pf}

We now verify  
% (specify the values of $ \g, \e_0, N_0 $) 
% By Theorem \ref{cor1} it remains to prove
that $ {\cal C}_\e   $ has asymptotically full measure, i.e. (\ref{Cmeas}) holds, choosing
\be\label{N0ge}
\g := \e^{\a} \quad {\rm with} \quad \a := 1 / (S + 1) \, , \quad N_0 := 4 \g^{-1} \, ,
\ee
so that (\ref{smallsto}) is fulfilled for $ \e $ small enough.

The complementary set of $ {\cal C}_\e $ in $ \Lambda $ has measure
\begin{eqnarray}
|{\cal C}_\e^c| & \stackrel{ \eqref{Cinfty}, \eqref{Gscavo}}= & 
\Big| \bigcup_{k \geq 1} {\mathtt G}_{N_k}^c (u_{k-1}) 
\bigcup_{k \geq  1 }  ({\cal G}_{N_k}^0(u_{k-1}))^c   \bigcup \tilde {\cal G}^c    \bigcup
 {\bar {\cal G}}^c \Big| \nonumber \\
& \leq & \sum_{k \geq  1 } | {\mathtt G}_{N_k}^c(u_{k-1}) | +
\sum_{k \geq  1 }  |({\cal G}_{N_k}^0(u_{k-1}))^c| +   |{\tilde {\cal G}}^c | + |{\bar {\cal G}}^c | \nonumber \\
& \stackrel{(\ref{measGN0}), (\ref{measBN0}), (\ref{tautau0}), \eqref{tildeGM}, (\ref{calGg})} \leq & 
C  \sum_{k \geq  1 } N_k^{-1}   
+  C \g  \leq C'  (N_0^{-1}  +  \g) 
\stackrel{(\ref{N0ge})} \leq C'' \e^{\a } \label{Cefull}
\end{eqnarray}
implying (\ref{Cmeas}). Finally \eqref{normd} follows by \eqref{uC1} and 
% From $(S2)_n$, the first step of section \ref{final} and , we derive
\begin{eqnarray*}
\| u(\e , \l)\|_{s_1} &\leq& \|u_0\|_{s_1} + \sum_{k=1}^{\infty} \|u_k -u_{k-1}\|_{s_1} \\
& \stackrel{(S1)_0, (S2)_n}\leq & N_0^{-\s} + \sum_{k=1}^{\infty}  N_k^{-\s-1} 
\leq C N_0^{-\s} \stackrel{\eqref{N0ge}}  \leq C \e^{\a \s} \, ,
\end{eqnarray*}
hence  $ \| u(\e,\l) \|_{s_1} \to 0 $, uniformly for $ \l \in \Lambda $,  as $ \e \to 0 $. 
Theorem \ref{thm:main}  is proved with $ s(d,\nu) := s_1 $
defined in \eqref{Sgr} and $ q(d,\nu) := S + 3 $, see \eqref{defk}. 
The $ C^\infty$-regularity result follows as in  \cite{BB10}-section 7.3.

\section{Appendix: proof of the Nash-Moser Theorem \ref{cor1}} \label{final}
\setcounter{equation}{0}

\noindent
{\bf Step 1: Initialization.} 
\label{sec:ini} We take $ \l \in {\cal N}(\bar {\cal G}, 2 N_0^{-\s} )$ (see \eqref{diofs}),
so that 
$$
{\cal L}_0 :=  P_0 (L_{\l \bar \om})_{| H_0}  \qquad {\rm satisfies} \qquad
\| {\cal L}_0^{-1} \|_{s_1} \leq  2N_0^{ \tau_1 + s_1} \g^{-1} 
$$
(see Lemma 7.1 in \cite{BB10}), and
we look for a solution of equation ($P_0$) as a fixed point of
$$
F_0 : H_0 \to H_0 \, , \quad F_0 (u) :=  \e {\cal L}_0^{-1} P_0 f(u) \, .  
$$
A contraction mapping argument (as in Lemma 7.2 of \cite{BB10}) proves that, 
for $\e \g^{-1} N_0^{\tau_1 +s_1+\s} \leq c(s_1) $ small, 
 $ \forall \l  \in  {\cal N}(\bar {\cal G}, 2N_0^{-\s}) $, 
there exists a unique solution $ {\wtilde u}_0 (\e, \l ) $ of $(P_0)$ in  
$$ 
{\mathtt B}_0(s_1) := \{ u \in H_0 \, : \, \| u \|_{s_1} \leq \rho_0 := N_0^{-\s} \} \, . 
$$
By uniqueness $ {\wtilde u}_0(0, \l ) = 0 $. 
The implicit function theorem implies that
$ {\wtilde u}_0(\e, \cdot ) \in C^1(  {\cal N}(\bar {\cal G}, 2N_0^{-\s}); H_0) $ and
$ \partial_{\l} {\wtilde u}_0 = -  {\cal L}_0^{-1}(\e)  (\partial_\l {\cal L}_0) {\wtilde u}_0 $
satisfies
$$
\| \partial_{\l} {\wtilde u}_0\|_{s_1} \leq 
C N_0^{\tau_1+ s_1+ 2 - \s}\g^{-1} \, .
$$
Then  we define the $ C^1 $ map $ u_0 := \psi_0 {\wtilde u}_0 : \Lambda \to H_0 $ with
cut-off function $ \psi_0 :  \Lambda \to [0,1] $,  
$$
\psi_0 :=
\begin{cases} 
1 \quad   {\rm if}   \  \l   \in  {\cal N}(\bar {\cal G}, N_0^{- \s }) \\
0  \quad  {\rm if} \  \l   \notin   {\cal N}(\bar {\cal G}, 2N_0^{- \s }) 
\end{cases} 
\quad  {\rm and} \qquad 
|D_\l \psi_0 | \leq N_0^\s C \, . 
$$
We get $ \| u_0 \|_{s_1} \leq N_0^{-\s} $, $ \| \partial_{\l} u_0 \|_{s_1} \leq C(s_1) N_0^{\t_1+ s_1 + 1} \g^{-1} $.
The statements $ (S1)_0 $,  $ (S4)_0 $ are proved (note that $ {\cal C}_0 = \tilde {\cal G}  \cap \bar {\cal G}  $). 
Statement $ (S5)_0 $ follows in the same way using \eqref{smallsto}.
Note that  $ (S2)_0 $, $ (S3)_0 $ are  empty. 

\smallskip 

For the next steps of the induction, the following lemma establishes a property which replaces $ (S3)_n $
for the first steps. It is proved exactly as in Lemma 7.3 of \cite{BB10}.

\begin{lemma}\label{Inizioind}
There exists $  N_0 (S, V) \in \N $ and  $ c(s_1) > 0 $ such that, if $N_0 \geq N_0(S,V)$ and  
$ \e N_0^{\t' + \d s_1}  \leq  c(s_1) $, 
then $ \forall N_0^{1/C_2} \leq N \leq N_0 $, $ \forall \| u \|_{s_1} \leq 1 $, we have $ {\cal G}_{N}(u) = \Lambda $. 
\end{lemma}

\smallskip

\noindent
{\bf Step 2: Iteration of the Nash-Moser scheme.} Suppose, by induction, that we have already defined 
$ u_n \in C^1(\Lambda; H_n ) $ 
and that properties $(S1)_k$-$(S5)_k $ hold for all $k\leq n$. 
We are going to define $u_{n+1}$ and prove the statements $(S1)_{n+1}$-$(S5)_{n+1}$.

\smallskip

In order to carry  out a modified Nash-Moser scheme, we shall study  the invertibility
of 
\be\label{caln+1}
{\cal L}_{n+1}(u_{n}) := P_{n+1} {\cal L}(u_{n})_{| H_{n+1}} \quad
{\rm where} \quad {\cal L}(u) := L_\om - \e (Df)(u) \, , 
\ee
(see  \eqref{Linve}) and the tame estimates of its inverse,   applying  Proposition \ref{propinv}. We  distinguish two cases. 
\\
If $ 2^{n+1} > C_2 $ (the constant $ C_2 $ is fixed in (\ref{tautau0})), 
then there exists a unique $ p \in [0,n] $ such that
\be\label{n+1np}
N_{n+1}  = N_p^{\chi}  \, , \quad  \chi = 2^{n+1-p} \in  [C_2, 2 C_2) \, , \  {\rm and } \quad 
N_{n+1} - 2 L_0  = N_p^{\tilde \chi}  \, , \quad  \tilde \chi  \in  [C_2, 2 C_2) \, . 
\ee
If $ 2^{n+1} \leq C_2 $ then there exists  $ \chi , \tilde \chi \in [C_2, 2 C_2] $ such that 
\be\label{Nbasso}
N_{n+1} = {\bar N}^\chi \, , \ \  {\bar N} := [ N_{n+1}^{1/C_2}] \in (N_0^{1/C_2}, N_0) \ \  {\rm and } \ \    
N_{n+1} - 2 L_0 = {\bar N}^{\tilde \chi} \, . %  \ \  {\bar N} := [ N_{n+1}^{1/C_2}] \in (N_0^{1/\chi}, N_0) \, . 
\ee
If \eqref{n+1np} holds we  consider  in Proposition \ref{propinv} the two scales $ N' = N_{n+1} $
(resp. $ N' = N_{n+1} - 2 L_0 $),
$ N = N_p $,  see (\ref{newscale}) .  %, resp. $ N' = N_{n+1} - 2 L_0 $, $ N = N_p $, and 
%$ \tilde \chi $. 
If \eqref{Nbasso}  holds, we set $N' = N_{n+1} $ (resp. $ N' = N_{n+1} - 2 L_0 $), $ N = \bar N $.

\begin{lemma}\label{H1H3} 
Let $ A(\e,\l,\teta) $ be defined in \eqref{matrpar}, with $u=u_n$. For all 
$$ 
\l \in \bigcap_{k =  1}^{n+1} {\cal G}_{N_k}^0 (u_{k-1}) \cap  \tilde {\cal G} 
 \, , \ %\forall 
\teta \in \R \, ,  
$$ 
the hypothesis  (H3) of Proposition \ref{propinv} apply
to $ A_{M,j_0}(\e,\l,\teta ) $,  $\forall M \in \{N_{n+1}, N_{n+1} - 2 L_0 \} $, 
$ \forall  j_0 \in \Z^d \setminus {\cal Q}_{M} $.  
\end{lemma}

\begin{pf}
We give the proof when $ M = N_{n+1}$ and  \eqref{n+1np} holds.  Since $ j_0 \notin {\cal Q}_{N_{n+1}} $
(i.e. $ (0,j_0) \notin \check {\cal Q}_{N_{n+1}} $) Lemma \ref{ANgood} implies that, a site 
\be\label{defE}
i \in E := (0,j_0) + [-N_{n+1}, N_{n+1}]^b
\ee
which is $ N_p $-good for 
$ A(\e, \l, \teta ) $ (see Definition \ref{GBsite}) %  with $ A = A(\e,\l,\teta) $)
 is  also  $ (A_{N_{n+1},j_0}(\e, \l, \teta), N_p)$-good 
(see Definition \ref{ANreg}).  
%with  $ A = A_{N_{n+1},j_0}(\e, \l ,\teta) $). 
As a consequence,
\be\label{badincl}
\Big\{ \ (A_{N_{n+1},j_0}(\e, \l, \teta), N_p){\rm -bad \ sites} \ \Big\}  \
\subset \  \Big\{ N_p{\rm -bad \ sites \ of} \  A(\e, \l ,\teta) \ {\rm with \ } |l| \leq N_{n+1} \Big\}.
\ee
and (H3) is proved if the latter $ N_p $-bad sites (in the right hand side of (\ref{badincl})) 
are contained in a disjoint union $ \cup_\a \Om_\a $ of clusters satisfying (\ref{sepabad}) (with $ N = N_p $). 
This is a consequence of  Proposition \ref{prop:separation} applied to 
the infinite dimensional matrix  $ A(\e, \l ,\teta) $. 
We claim that 
\be\label{inclup}
\bigcap_{k =  1}^{n+1} {\cal G}_{N_k}^0 (u_{k-1}) \subset {\cal G}_{N_p}(u_n) \, , \ {\rm i.e.} \
{\rm any} \ \, \l \, \in \bigcap_{k =  1}^{n+1} {\cal G}_{N_k}^0 (u_{k-1}) \ \, {\rm is} \ \, N_p-{\rm good \ for} \  \, A(\e,\l, \teta) \, ,  
\ee
and then assumption (i) of Proposition \ref{prop:separation} holds. Indeed, if $p=0$ then \eqref{inclup} 
is trivially true  because $ {\cal G}_{N_0}(u_n) =  \Lambda $, 
by Lemma \ref{Inizioind} and $(S1)_n $. If $ p \geq 1 $, we have
$$
\| u_n - u_{p-1} \|_{s_1} \leq \sum_{k=p}^{n} \| u_k - u_{k-1}  \|_{s_1} 
\stackrel{(S2)_k} \leq  \sum_{k=p}^{n} N_k^{-\s-1} \leq
N_p^{-\s}  \sum_{k\geq p} N_k^{-1} \leq N_p^{-\s} 
$$
and so $ (S3)_p $ implies 
$$
\bigcap_{k =  1}^p {\cal G}_{N_k}^0 (u_{k-1})  \subset {\cal G}_{N_p}(u_n) \, .
$$
Assumption (ii) of Proposition \ref{prop:separation} holds by (\ref{tautau0}), since $ \chi \in [C_2, 2 C_2) $.
Assumption (iii) of Proposition \ref{prop:separation} holds for all 
$ \l \in \tilde  {\cal G} $, see \eqref{ditilde}.

When \eqref{Nbasso} holds the proof is analogous using Lemma \ref{Inizioind} with $ N = \bar N $ and  $ (S1)_n $. 
\end{pf}

\begin{lemma}\label{S3n+1}
Property $ (S3)_{n+1} $ holds.
\end{lemma}

\begin{pf}
We want to prove that 
$$ 
\| u - u_n \|_{s_1} \leq N_{n+1}^{-\s} \ \ {\rm and} \ \  \l \in \bigcap_{k =  1}^{n+1} {\cal G}_{N_k}^0 (u_{k-1})
\cap {\tilde {\cal G} } 
\quad \Longrightarrow \quad  \l \in {\cal G}_{N_{n+1}} (u) \, .
$$ 
Since $  \l \in {\cal G}^0_{N_{n+1}} (u_n)$, 
by (\ref{BNcomponent2}) and Definition \ref{def:freqgood} 
it is sufficient to prove that 
%$ \forall M = N_{n+1}, N_{n+1} - 2 L_0 $,  $ \forall j_0 \in \Z^d \setminus {\cal Q}_{M} $,  
$$
B_{M} (j_0; \l)(u)  \subset B_{M}^0 (j_0;\l)(u_n)  \, , \quad 
\forall M \in \{ N_{n+1}, N_{n+1} - 2 L_0 \} \, ,  \ j_0 \in \Z^d \setminus {\cal Q}_{M} 
$$
(we highlight the dependence of these sets on $ u $, $ u_n $) 
or, equivalently, by (\ref{tetabadweak}), (\ref{tetabad}), that 
\be\label{inclusion+1}
(\| A_{M,j_0}^{-1} ( \e, \l , \teta )(u_n) \|_0 \leq M^\t  \  \Longrightarrow \
A_{M,j_0} ( \e, \l , \teta )(u) \ {\rm is} \ M - {\rm good}) \, ,  \ \forall  M \in \{ N_{n+1}, N_{n+1} - 2 L_0 \} \, ,
\ee
where $ A(\e, \l , \teta)(u ) $ is in \eqref{matrpar}. %  = {\cal L}(u) + \teta Y = L_\om  + \teta Y - \e (Df)(u) $.

Let us make the case $ M = N_{n+1} $, the other is similar. 
We prove (\ref{inclusion+1}) applying Proposition \ref{propinv} to $ A := A_{N_{n+1},j_0} ( \e, \l , \teta )(u)$
with $ E $ defined in \eqref{defE},  $  N' = N_{n+1} $, $ N =  N_p $ (resp. $N=\bar N$) if  \eqref{n+1np}
(resp. \eqref{Nbasso}) is satisfied. 

Using Lemma \ref{lem:multi}, $ \| V \|_{C^q} \leq C $,  
assumption (H1)  holds with 
\be\label{decadimentoA}
\Upsilon  \leq C (1+ \| u_n \|_{s_1}+ \norma V \norma_{s_1}) 
\stackrel{ (S1)_n} \leq C' (V) \, .  
\ee  
By Lemma \ref{H1H3}, for all $  \teta \in \R $, $ j_0 \in \Z^d \setminus {\cal Q}_{N_{n+1} } $, 
the hypothesis  (H3) of Proposition \ref{propinv} holds for $ A_{N_{n+1},j_0}(\e, \l ,\teta )(u_n) $. 
Hence, by Proposition \ref{propinv}, for $s\in [s_0,s_1]$, if
$$
\| A_{N_{n+1},j_0}^{-1} ( \e, \l , \teta )(u_n) \|_0 \leq N_{n+1}^\t  
$$
(which is assumption (H2)) then
\be\label{risultAN}
\nors{A^{-1}_{N_{n+1},j_0}(\e, \l ,\teta )(u_n)} \leq 
\frac{1}{4} N_{n+1}^{\tau' } \Big( N_{n+1}^{\d s}+  \norma V \norma_s + \e \nors{ (Df)(u_n) } \Big) \, .
\ee
Finally, since $ \| u - u_n \|_{s_1} \leq N_{n+1}^{-\sigma} $ we have
$$
\norsone{A_{N_{n+1},j_0} ( \e, \l , \teta )(u_n)-A_{N_{n+1},j_0} ( \e, \l , \teta )(u)} \leq C\e \| u - u_n \|_{s_1} \leq N_{n+1}^{-\s} 
$$
and  \eqref{inclusion+1} follows by \eqref{risultAN} and 
a standard perturbative argument (see e.g. \cite{BB10}). 
%(see for instance \eqref{inv1} in  Lemma \ref{leftinv} with any $s\in [s_0,s_1]$ instead of $s_0$). 
\end{pf}

From now on the convergence proof of the Nash-Moser iteration follows \cite{BB10} with no changes.

In order to define $ u_{n+1} $, we write, for $ h \in H_{n+1} $, 
\be\label{scritt1}
P_{n+1}  \Big(L_\om ( u_n + h ) - \e f( u_n + h ) \Big) 
=  r_n + {\cal L}_{n+1} (u_n) h + R_n ( h )  
\ee
where 
$ {\cal L}_{n+1}(u_{n}) $
is defined in  \eqref{caln+1} and 
\be\label{Rnh}
r_n :=  P_{n+1} \Big(L_\om u_n - \e  f( u_n ) \Big) \, ,  \quad 
R_n ( h ) := - \e  P_{n+1} \Big( f( u_n + h ) - f( u_n ) - (Df)( u_n ) h \Big) \, .
\ee
By $ (S4)_n $, if $ \l \in {\cal N}({\cal C}_n, N_{n}^{-\s}) $ then $ u_n $ solves Equation $(P_n) $ and so
\be\label{rnhigh}
r_n =  P_{n+1}P_n^\bot  \Big( L_\om  u_n - \e   f( u_n )\Big)  = 
P_{n+1}P_n^\bot  \Big(V_0 \,  u_n - \e  f( u_n ) \Big) \, , 
\ee
using also that $ P_{n+1} P_n^\bot ( D_\om u_n) = 0 $, see \eqref{Lomega}.
Note that, by (\ref{defNn}) and $ \s \geq 2 $ (see \eqref{def:sigma}), for $ N_0  \geq 2 $, we have the inclusion
$ {\cal N}({\cal C}_{n+1}, 2 N_{n+1}^{-\s}) \subset {\cal N}({\cal C}_n, N_{n}^{-\s})$.

\begin{lemma}\label{invLn+1}
{\bf (Invertibility of $ {\cal L}_{n+1} $)}  
For all $  \l  \in  {\cal N}({\cal C}_{n+1}, 2 N_{n+1}^{-\s}) $
the operator  $ {\cal L}_{n+1}(u_n) $ is invertible and,  for $ s  = s_1, S $, 
\be\label{goal}
\norma {\cal L}_{n+1}^{-1}(u_n) \norma_s \leq  N_{n+1}^{\t' + \d s } \,  . 
\ee
As a consequence, by (\ref{opernorm}),  $ \forall h \in H_{n+1} $, 
\be\label{os1}
\| {\cal L}_{n+1}^{-1}(u_n) h \|_{s_1} \leq C(s_1) N_{n+1}^{\t' + \d s_1 } \| h \|_{s_1} \, , 
\ee
\be\label{os2}
\| {\cal L}_{n+1}^{-1}(u_n) h \|_S \leq  N_{n+1}^{\t' + \d s_1 } \|h\|_S + C(S) N_{n+1}^{\t' + \d S }  \|h \|_{s_1} \, .
\ee
\end{lemma}

\begin{pf} 
We  apply the multiscale Proposition \ref{propinv} to $ A_{N_{n+1}} = {\cal L}_{n+1} (u_n) $ as in % as in \cite{BB10}, see
Lemma \ref{S3n+1}. Assumption (H1) holds by \eqref{decadimentoA}. For all
 $ \l \in {\mathtt G}_{N_{n+1}}(u_n) $ (see (\ref{Binver})) 
$ \| {\cal L}_{n+1}^{-1}(u_n) \|_0 \leq N_{n+1}^\tau $ and (H2) holds.
The hypothesis (H3)   holds, for $ \l \in {\cal C}_{n+1}$ (see \eqref{Gscavo}), 
as a particular case of Lemma \ref{H1H3}, for $ \teta = 0 $,  $ j_0 = 0 $, $ M = N_{n+1}$, 
and since $ 0 \notin {\cal Q}_{N_{n+1}} $.
Then Proposition \ref{propinv} applies $ \forall \l \in {\cal C}_{n+1} $, implying \eqref{goal}.
 For all $ \l \in {\cal N}({\cal C}_{n+1}, 2 N_{n+1}^{-\s}) $ the proof of \eqref{goal} 
follows by a perturbative argument  as in Lemma 7.7 in \cite{BB10}.
\end{pf}

By \eqref{scritt1}, the equation ($ P_{n+1} $) is equivalent to the fixed point problem $ h = F_{n+1}( h ) $ 
where
$$
F_{n+1}: H_{n+1} \to H_{n+1} \, ,  \qquad
F_{n+1}(h) := - {\cal L}_{n+1}^{-1} (u_n) ( r_n +  R_n ( h )) \, .
$$
% \begin{lemma} \label{lemcon}
By a contraction mapping argument as in Lemma 7.8 in \cite{BB10}
(using \eqref{os1}, % on the inverse $ {\cal L}_{n+1}^{-1} (u_n) $ and 
 \eqref{rnhigh}, \eqref{Rnh})
we prove the existence,  
$ \forall  \l \in {\cal N}({\cal C}_{n+1}, 2 N_{n+1}^{-\s}) $, of a unique fixed point $ {\wtilde h}_{n+1}(\e, \l) $
of $ F_{n+1} $ 
in
$$
{\mathtt B}_{n+1}(s_1) := \Big\{ h \in  H_{n+1} \  : \  \|h\|_{s_1} \leq  \rho_{n+1} :=    N_{n+1}^{-\s-1} \Big\} \, . 
$$
Since $ u_n (0,\l) = 0 $ (by $ (S1)_n $), we deduce, by the uniqueness of the fixed point, that
$ {\wtilde h}_{n+1}(0,\l) = 0 $.
% \quad \forall (0,\l) \in {\cal N}({\cal C}_{n+1}, 2N_{n+1}^{-\s}) $. 
Moreover, as in Lemma 7.9 of \cite{BB10} (using the tame estimate \eqref{os2}),
one deduces  the following bound on the high norm 
$$ 
\|  {\wtilde h}_{n+1}  \|_{S} \leq K(S) N_{n+1}^{\t' + \d s_1} U_n \, . 
$$
By the implicit function theorem as in Lemma 7.10 in \cite{BB10}
(using \eqref{os1}-\eqref{os2}) 
the map $ {\wtilde h}_{n+1}  $ is in $  C^1({\cal N}({\cal C}_{n+1}, 2 N_{n+1}^{-\s}), H_{n+1}) $ and 
$$
\| \partial_{\l} {\wtilde h}_{n+1} \|_{s_1} \leq  N_{n+1}^{ -1}  \, , \quad 
 \| \partial_{\l} {\wtilde h}_{n+1} \|_{S} \leq  N_{n+1}^{\t'+\d s_1 +1} 
 \Big(N_{n+1}^{\t'+\d s_1 +1} U_n + U_n' \Big) \, . 
$$
Finally we define 
the $ C^1 $-extension 
% of $ ({\wtilde h}_{n+1})_{|{\cal C}_{n+1}}  $   
onto the whole $  \Lambda $ as
$$
h_{n+1} (\l) := \left\{ \begin{array}{lll} 
\psi_{n+1}(\l) {\wtilde h}_{n+1}(\l) & {\rm if} &  \l \in {\cal N}({\cal C}_{n+1}, 2 N_{n+1}^{- \s } )   \\
0 & {\rm if} & \l \notin  {\cal N}({\cal C}_{n+1},2 N_{n+1}^{-\s} )  \end{array}
\right.
$$
where $ \psi_{n+1} $ is a $C^\infty $ cut-off function satisfying 
$$ 
0 \leq \psi_{n+1} \leq 1 \, , \quad
\psi_{n+1} \equiv 
\begin{cases} 
1 \quad {\rm if} \  \l \in {\cal N}({\cal C}_{n+1},  N_{n+1}^{-\s} ) \\
0  \quad  {\rm if} \  \l \notin  {\cal N} ({\cal C}_{n+1}, 2{N_{n+1}^{- \s})} 
\end{cases}  {\rm and} \quad 
|\partial_{\l} \psi_{n+1} | \leq N_{n+1}^\s C \, . 
$$
Then  %satisfies 
(see Lemma 7.11 in \cite{BB10})
$$ 
\| h_{n+1}\|_{s_1} 	\leq N_{n+1}^{-\s-1} \, , \quad 
\| \partial_{\l} h_{n+1}\|_{s_1}	 \leq   N_{n+1}^{-1/2}  \, .  
$$
% and $ h_{n+1} $ is equal to $ {\wtilde h}_{n+1} $ on ${\cal N}({\cal C}_{n+1}, N_{n+1}^{-\s}  ) $. 
In conclusion,  
$ u_{n+1} :=  u_n + h_{n+1} $
satisfies
$ (S1)_{n+1} $, $ (S2)_{n+1} $, $ (S4)_{n+1} $, $ (S5)_{n+1} $ (see Lemma 7.12 in \cite{BB10}).

\noindent
Massimiliano Berti, Dipartimento di Matematica e Applicazioni ``R. Caccioppoli",
Universit\`a degli Studi Napoli Federico II,  Via Cintia, Monte S. Angelo, 
I-80126, Napoli, Italy,  {\tt m.berti@unina.it}.
\\[2mm]
Philippe  Bolle, Universit\'e
d'Avignon et des Pays de Vaucluse, Laboratoire d'Analyse non Lin\'eaire 
et G\'eom\'etrie (EA 2151), F-84018 Avignon, France, {\tt philippe.bolle@univ-avignon.fr}.
\\[2mm]
\indent
This research was supported by the European Research Council under FP7 ``{\sl New Connections between Dynamical Systems
and Hamiltonian PDEs with Small Divisors Phenomena}" and partially by
the PRIN2009 grant ``{\sl Critical Point Theory and Perturbative Methods for Nonlinear Differential Equations}".

\begin{thebibliography}{10}


\bibitem{BBi} Berti M., Biasco L., {\it Branching of Cantor manifolds of elliptic tori
and applications to  PDEs}, Comm. Math. Phys, 305, 3,  741-796, 2011.

\bibitem{BBP1} Berti M., Biasco L., Procesi M., {\it KAM theory for
Hamiltonian derivative wave equations}, preprint 2011.

\bibitem{BB07} Berti M., Bolle P., 
{\it Cantor families of periodic solutions of wave equations with $ C^k $ nonlinearities},  NoDEA
Nonlinear Differential Equations Appl., 15, 247-276, 2008.

\bibitem{BB}
Berti M., Bolle P., {\it Sobolev Periodic solutions of nonlinear wave equations in higher spatial dimension}, 
Archive for Rational Mechanics and Analysis, 195, 609-642, 2010.

\bibitem{BB10}
Berti M., Bolle P., {\it Quasi-periodic solutions  with Sobolev regularity of NLS on $ \T^d $ with a multiplicative potential}, 
to appear on Journal European Math. Society.

\bibitem{BBP}
Berti M., Bolle P., Procesi M., {\it An abstract Nash-Moser theorem with parameters 
and applications to PDEs}, Ann. I. H. Poincar\'e, 27, 377-399, 2010.


\bibitem{BP}
Berti M., Procesi M., {\it Nonlinear wave and Schr\"odinger equations on compact Lie groups and homogeneous spaces}, 
Duke Math. J., 159, 3, 479-538, 2011.


\bibitem{Bo1}  Bourgain J., {\it  Construction of quasi-periodic solutions 
for Hamiltonian perturbations of linear equations and applications 
to nonlinear PDE}, Internat. Math. Res. Notices, no. 11, 1994.


\bibitem{B4} Bourgain J., {\it Construction of periodic solutions of nonlinear
wave equations in higher dimension},  Geom. Funct. Anal. 5,  no. 4, 629-639, 1995.

\bibitem{B3}  Bourgain J., {\it Quasi-periodic solutions of Hamiltonian
perturbations of $2D$ linear Schr\"odinger equations},
Annals of Math. 148, 363-439, 1998.

\bibitem{B5}  Bourgain J.,  {\it Green's function estimates for lattice Schr\"odinger 
operators and applications}, Annals of Mathematics Studies 158, 
Princeton University Press, Princeton, 2005.

\bibitem{BW1} Bourgain J., Wang W.M.,  
{\it Anderson localization for time quasi-periodic
random Schr\"odinger and wave equations}, Comm. Math. Phys. 248, 429 - 466, 2004.


\bibitem{C} Craig W.,  {\it Probl\`emes de petits diviseurs dans
les \'equations aux d\'eriv\'ees partielles},
Panoramas et Synth\`eses, 9, 
Soci\'et\'e Math\'ematique de France, Paris, 2000. 

\bibitem{CW} Craig W., Wayne C. E., {\it Newton's method and periodic solutions
of nonlinear wave equation}, Comm. Pure  Appl. Math. 46, 1409-1498, 1993.

\bibitem{EK} Eliasson L. H., Kuksin S.,  {\it   KAM for nonlinear Schr\"odinger equation}, 
Annals of Math., 172, 371-435, 2010.

\bibitem{EK1} Eliasson L. H., Kuksin S., {\it On reducibility of Schr\"odinger equations with
quasiperiodic in time potentials}, Comm. Math. Phys, 286, 125-135, 2009.

%\bibitem{Fo} Fokam J. M. , 
%{\it Forced vibrations via Nash-Moser iteration}, Comm. Math. Phys. 283, 2, 285-304, 2008.

\bibitem{GXY} Geng J., Xu X.,  You J., {\it An infinite dimensional KAM theorem and 
its application to the two dimensional cubic 
Schr\"odinger equation}, Adv. Math. 226, 6, 5361-5402, 2011.

\bibitem{K1} Kuksin S., {\it Hamiltonian perturbations
of infinite-dimensional linear systems with imaginary spectrum},
Funktsional Anal. i Prilozhen. 2, 22-37, 95, 1987.

\bibitem{K2} Kuksin S., {\it Analysis of Hamiltonian PDEs}, Oxford
Lecture series in Mathematics and its applications 19, Oxford University
Press, 2000.

\bibitem{KP} Kuksin S., P\"oschel J., {\it
Invariant Cantor manifolds of quasi-periodic oscillations
for a nonlinear Schr\"{o}dinger equation}, Annals of Math. (2)  143, 149-179, 1996.


%\bibitem{PY}  Plotnikov P.I., Yungerman  L.N.,
%\emph{Periodic solutions of a weakly nonlinear wave equation
%with an irrational relation of period to interval length}, 
%(Russian) Differentsial nye Uravneniya 24 (1988), 
%no. 9, 1599--1607, 1654;
%translation in Diff. Equations \textbf{24} (1988), no. 9, 
%1059--1065 (1989).

\bibitem{Po2} P\"oschel J.,
{\it A KAM-Theorem for some nonlinear partial differential equations},  Ann. Scuola Norm. Sup.%
Pisa Cl. Sci.(4), 23, 119-148, 1996.

\bibitem{Po3} P\"oschel J., {\it Quasi-periodic solutions for
a nonlinear wave equation}, Comment. Math. Helv.,  71,
no. 2, 269-296, 1996.

\bibitem{PX} Procesi M., Xu X., {\it Quasi-T\"oplitz Functions in KAM Theorem}, preprint 2011.


\bibitem{W1} Wang W. M.,
{\it Supercritical nonlinear Schr\"odinger equations I: quasi-periodic solutions},  preprint 2010.

\bibitem{W2} Wang W. M.,
{\it Supercritical nonlinear wave equations: quasi-periodic solutions and almost global existence}, 
preprint 2011.

\bibitem{Wa1}  Wayne E., {\it Periodic and quasi-periodic solutions
of nonlinear wave equations via KAM theory},
Comm. Math. Phys. 127, 479-528, 1990.


\end{thebibliography}
\end{document}